\documentclass[10pt,twoside,reqno]{amsart} 
\usepackage{amsmath, amsfonts, amsthm, amssymb}
\usepackage[francais,english]{babel}
\usepackage[latin1]{inputenc}
\usepackage{amsmath}
\usepackage{amsfonts}
\usepackage{amsthm}
\usepackage[]{graphics}
\usepackage{color}
\usepackage[T1]{fontenc}
 \usepackage{setspace}
\usepackage{algorithm}
\usepackage{algorithmic}

\usepackage[T1]{fontenc}
\usepackage[latin1]{inputenc}
\usepackage{amsmath,latexsym}
\usepackage{here}
\usepackage{graphicx}
\usepackage{multirow}
\usepackage{amssymb}

\definecolor{red}{rgb}{0.82,0.00,0.00}  

\numberwithin{equation}{section}
\input epsf
 

\newcommand{\ds}{\displaystyle}
\def\R{{\rm I\hspace{-0.50ex}R} }
\def\P{{\rm I\hspace{-0.50ex}P} }

\def\w{{\bf w}}
\def\f{{\bf f}}
\def\v{{ \bf v}}
\def\u{{\bf u}}
\def\b{{\bf b}}
\def\0{{\bf 0}}

\def\x{{\bf{x}}}

\def\n{{\bf n}}
\def\div{\operatorname{div}}

\newtheorem{lem}{Lemma}[section]
\newtheorem{conc}[lem]{Conclusion}
\newtheorem{thm}[lem]{Theorem}
\newtheorem{rmq}[lem]{Remark}

\renewcommand{\epsilon}{\varepsilon}
\def\twoplot[#1]#2#3#4#5{
\begin{figure}[hbt]
\begin{multicols}{2}
\begin{center}
    \includegraphics*[#1]{#2}
    \caption{\label{#2} #4}
\end{center}
\begin{center}
    \includegraphics*[#1]{#3}
    \caption{\label{#3} #5}
\end{center}
\end{multicols}
\end{figure}
}

\begin{document}

\bibliographystyle{plain}

\title[Convergence analysis of numerical schemes]{ Convergence analysis of numerical schemes for the Darcy-Forchheimer problem}
\author[ SAYAH  ]{Toni Sayah$^{\dagger}$}
\thanks{ \today.
\newline
$^{\dagger}$  Laboratoire de "Math\'ematiques et applications", Unit\'e de recherche "Math\'ematoqies et Mod\'elisation", CAR, Facult\'e des sciences, Universit\'e Saint-Joseph, Lebanon.
\newline
toni.sayah@usj.edu.lb.}


\begin{abstract}
This paper deals with the Darcy-Forchheimer problem with two kinds of boundary conditions. We discretize the system by using the finite element methods and we propose two iterative schemes to solve the discrete problems. The well-posedness and the convergence of the corresponding iterative problems are then proven. Finally,  several numerical experiments are performed to validate the proposed numerical schemes.
\end{abstract}

\maketitle

\section{Introduction}\label{intro}
Darcy's law describes the creeping flow of Newtonian fluids in porous media. It gives the following linear relationship between the velocity of creep flow and the gradient of pressure
\[
\ds \frac{\mu}{\rho} K^{-1}  \u + \nabla p = \f.
\]
{  A theoretical derivation of Darcy's law can be found in \cite{Neum,Whit}.} Forchheimer \cite{Forchheimer} conducted flow experiments in sandpacks and recognized that when the velocity is higher and the porosity is nonuniform, Darcy's law is not adequate. He proposed the following nonlinear equation:
\[
\ds \frac{\mu}{\rho} K^{-1}  \u + \frac{\beta}{\rho} |\u| \u+ \nabla p = \f.
\]
A theoretical derivation of Forchheimer's law can be found in \cite{RM}.\\
\noindent A mixed element for Forchheimer equation (or called Darcy-Forchheimer equation sometimes) was introduced by Girault and Wheeler \cite{GW}. They proved the existence and uniqueness of a weak solution for the Forchheimer equation. At the discrete level, they approximated velocity by piecewise constants and pressure by Crouzeix-Raviart element. They also proposed an alternating directions iterative method to solve the system of nonlinear equations obtained by finite element discretizaton. The convergence of both the iterative algorithm and the mixed element scheme are presented, and the error estimate of the mixed element scheme is demonstrated too. Lopez et al. \cite{LMJ} carried out numerical tests of the methods studied in \cite{GW} in order to corroborate the results presented there. Furthermore, another mixed finite element space was proposed in which the approximation to the pressure is smoother than the one obtained with the space proposed in \cite{GW}. In \cite{HH}, the authors proposed a mixed element approximation: the Raviart-Thomas mixed element, Brezzi-Douglas-Marini mixed element. They demonstrated the existence and uniqueness of the weak solution, and gave the corresponding error estimate. They also introduced an iterative algorithm, without the study of the corresponding convergence, and showed numerical tests. Salas J. et al. \cite{SJ} presented a theoretical study of the mixed finite element space, as proposed in \cite{LMJ}, and showed the discrete solution's existence and uniqueness, its convergence and an error estimate.

\noindent Let $\Omega$ be a bounded subset of $\R^d$ ($d = 2, 3$) with Lipschitz continuous boundary $\Gamma=\partial \Omega$. {  We} consider the Darcy-Forchheimer equation
\begin{equation}\label{E1}
\ds \frac{\mu}{\rho} K^{-1}  \u + \frac{\beta}{\rho} |\u| \u+ \nabla p = \f \;\;\;\;\; \mbox{in}\;\;\; \Omega,
\end{equation}
with the divergence constraint
\begin{equation}\label{E2}
\ds \div \u = b \;\;\;\;\; \mbox{in}\;\;\; \Omega,
\end{equation}
and the boundary condition
\begin{equation}\label{E3}
\ds \u \cdot \n = g_u \;\;\;\;\; \mbox{on}\;\; \partial \Omega ,
\end{equation}
or
\begin{equation}\label{E4}
\ds p = g_p \;\;\;\;\; \mbox{on}\;\; \partial \Omega.
\end{equation}
Here $\u$ represents the velocity, $p$ represents the pressure, $\n$ is the unit exterior normal vector to $\Gamma$, $|.|$ denotes the Euclidean norm, $|\u|^2 = \u \cdot \u$. {  The parameters} $\rho$, $\mu$ and $\beta$ represent the density of the fluid, its viscosity and its dynamic viscosity, respectively. $\beta$ is also referred as Forchheimer number when it is a scalar positive constant. $K$ is the permeability tensor, assumed to be uniformly positive definite and bounded such that there exists two positive real numbers $K_m$ and $K_M$ such that
\begin{equation}\label{KmM}
0 < K_m \, \x \cdot \x \le (K^{-1}(\x)\x) \cdot \x \le K_M \, \x \cdot \x.
\end{equation}
{  It should be noted that $K_m$ could be very close to zero and $K_M$ could be very large.}\\

\noindent We denote by Problem $(P_1)$ the system of equations (\eqref{E1}, \eqref{E2}, \eqref{E3}), and by Problem $(P_2)$ the system of equations (\eqref{E1}, \eqref{E2}, \eqref{E4}). For the compatibility, we suppose that $b$ and $g_u$ verify the following compatibility condition:
\[
\ds \int_\Omega b(\x) d\x = \ds \int_\Gamma \mbox{g}_u(s) \, ds.
\]
In this paper we consider Problems $(P_1)$ and $(P_2)$  and recall corresponding discrete schemes introduced in \cite{SJ} and \cite{HH}. {  Two discrete iterative schemes are proved and their corresponding convergences are showed}. The technique of the convergence is inspired from \cite{sayahjsc}. We also introduce numerical results of validation.
This paper is organised as follow :
\begin{itemize}
\item Section 2 describes the problem and the weak formulations.
\item Section 3 is devoted to study of the discretisation and the convergence of the proposed iterative schemes.
\item Section 4 is devoted to the numerical results.
\end{itemize}
\section{Notations and weak formulations}
In order to introduce the variational formulations, we  recall  some classical Sobolev spaces and  their  properties.

Let $\alpha=(\alpha_1,\alpha_2, \dots \alpha_d)$ be a $d$-uple of non negative
integers, set $|\alpha|=\ds \sum_{i=1}^d \alpha_i$, and define the partial
derivative $\partial^\alpha$ by
$$
\partial^{\alpha}=\ds \frac{\partial^{|\alpha|}}{\partial x_1^{\alpha_1}\partial
x_2^{\alpha_2}\dots\partial x_d^{\alpha_d}}.
$$
Then, for any positive integer $m$ and number {  $q\geq 1$}, we recall
the classical Sobolev space
\begin{equation}
\label{eq:Wm,p}
{ W^{m,q}}(\Omega)=\{v \in {  L^q}(\Omega);\,\forall\,|\alpha|\leq
m,\;\partial^{\alpha} v \in { L^q}(\Omega)\},
\end{equation}
equipped with the seminorm
\begin{equation}
\label{eq:semnormWm,p}
|v|_{{ W^{m,q}}(\Omega)}=\big\{\sum_{|\alpha|=m} \int_{\Omega}
|\partial^{\alpha} v|^{  q} \,d\x\,\big\}^{ \frac{1}{q}}
\end{equation}
and the norm
\begin{equation}
\label{eq:normWm,p}
\|v\|_{  W^{m,q}(\Omega)}=\big\{\sum_{0\leq k\leq m}
|v|_{  W^{k,q}(\Omega)}^{  q} \big\}^{  \frac{1}{q}}.
\end{equation}
When $  q=2$, this space is the Hilbert space $H^m(\Omega)$.
In
particular, the scalar product of $L^2(\Omega)$ is denoted by
$(.,.)$. Furthermore, we recall the following standard spaces for Darcy's equations:
\[
L^2_0(\Omega) = \{ v\in L^2(\Omega); \ds \int_\Omega v d\x =0 \}.
\]
The definitions of these spaces are extended straightforwardly to
vectors, with the same notation, but with the following
modification for the norms in the non-Hilbert case. Let $\v$ be a
vector valued function {  and we define the norm}
\begin{equation}
\label{eq:normLp} \|\v\|_{  L^q(\Omega)}= \big(\int _{\Omega}
|\v|^{  q}\,d\x\, \big)^{  \frac{1}{q}},
\end{equation}
where $|.|$ denotes the Euclidean vector norm.\\

We recall the following standard space
\begin{equation}
\label{eq:Hdiv}
H(\div,\Omega)=\{\v\in L^2(\Omega)^d;\,\div(\v)\in L^2(\Omega)\},
\end{equation}
equipped with the norm
\begin{equation}
\label{eq:normHdiv}
\|\v\|_{H(\div,\Omega)}^2=\|\v\|_{L^2(\Omega)}^2+\|\div(\v)\|_{L^2(\Omega)}^2.
\end{equation}
\noindent Let us now introduce the following technical lemma:
{\lem \label{relatinf} For all $x,y \in \R$ and {  $q\in \R^{+}$, {  the following bound holds:}
\[
(|x|^{q} x - |y|^{q} y)(x-y) \ge 0.
\]
}
}
%
%
\subsection{Variational formulation of the first problem}
We refer to \cite{GW} for all the properties and the details of the weak formulation corresponding to Problem $(P_1)$ presented in this section. {  Let} us introduce the spaces:
\[
\begin{array}{ll}
\medskip
X_u = L^3(\Omega)^d, \\
M_u = W^{1, \frac{3}{2}}(\Omega) \cap L^2_0(\Omega),
\end{array}
\]
which satisfy the following inf-sup condition {  (see \cite{GW} for the proof)}
\[
\underset{q\in M_u}{\inf} \; \; \underset{\v \in X_u}{\sup} \; \; \ds \frac{\ds \int_\Omega \nabla q \cdot \v \, d\x}{||\nabla q ||_{L^{3/2}(\Omega)} ||\v||_{L^3(\Omega)}} = 1.
\]
The velocity $\u$ and pressure $p$ of Problem $(P_1)$ are considered respectively in $X_u$ and $M_u$, therefore $b$ and $g_u$ are assumed to be respectively {  in $L^{3d/(d+3)}(\Omega)$ and $L^{3(d-1)/d}(\Gamma)$ }(see \cite{GW} for details). \\

Problem $(P_1)$ is equivalent to the following variational formulation: Find $(\u,p)\in X_u \times M_u$ such that
\begin{equation}\label{V1}
\left\{
\begin{array}{ll}
\medskip
\forall \v \in X_u, \quad \ds \frac{\mu}{\rho} \int_\Omega K^{-1} \u \cdot \v \, d\x + \frac{\beta}{\rho} \int_\Omega |\u| \u \cdot \v \, d\x + \int_\Omega \nabla p \cdot \v \, d\x = \ds \int_\Omega \f \cdot \v \, d\x,\\
\forall q\in M_u, \quad \ds \int_\Omega \nabla q \cdot \u \, d\x = -\int_\Omega b \, q \, d\x + \int_\Gamma g_u \, q \, ds.
\end{array}
\right.
\end{equation}
It can be demonstrated that, for all $b \in L^{3d/(d+3)}(\Omega)$ and {  $g_u \in L^{3(d-1)/d}(\Gamma)$}, there is only one $\u_l \in L^3(\Omega)^d/V_u$ such that
\[
\forall q \in M_u, \, \ds \int_\Omega \nabla q \cdot \u_l \, d\x = -\int_\Omega b \, q \, d\x + \int_\Gamma g_u \, q \, ds,
\]
and
\[
{  \ds ||\u_l ||_{L^3(\Omega)^d/V_u} \le C \big( ||b||_{L^{3d/(d+3)}(\Omega)}  + ||g_u ||_{L^{3(d-1)/d}(\Gamma)}  \big),}
\]
where
\[
V_u = \{ \v\in X_u; \forall q \in M_u, \ds \int_\Omega \nabla q \cdot \v \, d\x= 0  \},
\]
and $C$ is a constant depending only on $\Omega$ and $d$.\\
We have from \cite{GW} that Problem $(P_1)$ is equivalent to \eqref{V1} which admits a unique solution in $(\u,p) \in X_u \times M_u$ satisfying the relations
{
\begin{equation}
\begin{array}{ll}\label{bound1u}
\medskip
||\u||_{L^3(\Omega)} \le C \Big(  ||\u_l||_{L^3(\Omega)}^3 + ||\u_l||_{L^2(\Omega)}^2 + ||\f ||_{L^2(\Omega)}^2 \Big)^{1/3}, \\
||\nabla p||_{L^{3/2}(\Omega)} \le C \Big(  ||\u||_{L^{3/2}(\Omega)} +  ||\u ||_{L^3(\Omega)}^2 + ||\f ||_{L^2(\Omega)} \Big).
\end{array}
\end{equation}
}
\subsection{Variational formulation of the second problem}
In this section, we consider $g_p=0$ for the simplicity and refer to \cite{HH} for all the properties and the details of the weak formulation corresponding to Problem $(P_2)$ presented in this section. {  Let us introduce} the spaces $M_p = L^2(\Omega)$ and
\[
X_p = \{ \v \in L^3(\Omega)^d; \div (\v) \in L^2(\Omega)\},
\]
equipped with the norm
\[
||\v ||_{X_p} = ||\v ||_{L^3(\Omega)} + ||\div(\v)||_{L^2(\Omega)}.
\]
The spaces $X_p$ and $M_p$ {  satisfy} the following inf-sup condition: there exists a positive constant $\gamma$ such that
\[
\underset{q\in M_p}{\inf} \; \;  \underset{\v \in X_p}{\sup} \; \; \ds \frac{\ds \int_\Omega q \, \div(\v) \, d\x}{||q ||_{M_p} ||\v||_{X_p}} \ge \gamma.
\]
In this case, we assume $\f = \nabla Z \in L^2(\Omega)^d$ the gradient of the depth function $Z\in H^1(\Omega)$ and $b\in L^2(\Omega) $ (see \cite{HH} for details). \\

Problem $(P_2)$ is equivalent to the following variational formulation: Find $(\u,p)\in X_p \times M_p$ such that
\begin{equation}\label{V2}
\left\{
\begin{array}{ll}
\medskip
\forall \v \in X_p, \quad \ds \frac{\mu}{\rho} \int_\Omega K^{-1} \u \cdot \v \, d\x + \frac{\beta}{\rho} \int_\Omega |\u| \u \cdot \v \, d\x - \int_\Omega p \, \div(\v) \, d\x =  \int_\Omega \f \cdot \v \, d\x,\\
\forall q\in M_p, \quad \ds \int_\Omega  q \, \div(\u) \, d\x = \int_\Omega b \, q \, d\x .
\end{array}
\right.
\end{equation}
We have from \cite{HH} that Problem $(P_2)$ is equivalent to \eqref{V2} which admits a unique solution in $(\u,p) \in X_p \times M_p$ satisfying the relations
\begin{equation}\label{bound2u}
\begin{array}{ll}
\medskip
||\u||^2_{L^2(\Omega)} + ||\u||^3_{L^3(\Omega)} + ||\div(\u)||^2_{L^2(\Omega)} \le C(  ||b||^2_{L^2(\Omega)} + ||b||^3_{L^2(\Omega)} + ||\f||^2_{L^2(\Omega)}) , \\
|| p||_{L^2(\Omega)} \le {  C} ( \ds  ||b||_{L^2(\Omega)} + ||b||^2_{L^2(\Omega)} + ||\f||_{L^2(\Omega)} + ||\f||^2_{L^2(\Omega)}  ).
\end{array}
\end{equation}
\section{Finite element discretization and convergence}
From now on, we assume that $\Omega$ is a polygon when $d=2$ or
polyhedron when $d=3$, so it can be completely meshed.  For
 the space discretization,  we consider a
regular (see Ciarlet~\cite{PGC}) family of triangulations
$( \mathcal{T}_h )_h$ of $\Omega$ which is a set of closed non degenerate
triangles for $d=2$ or tetrahedra for $d=3$, called elements,
satisfying,
\begin{itemize}
\item for each $h$, $\bar{\Omega}$ is the union of all elements of
 $\mathcal{T}_h$;
\item the intersection of two distinct elements of  $\mathcal{T}_h$ is
either empty, a common vertex, or an entire common edge (or face
when $d=3$);
\item the ratio of the diameter   {  $h_{\kappa}$  of an element $\kappa \in \mathcal{T}_h$  to
the diameter  $\rho_\kappa$}  of its inscribed circle when $d=2$ or ball when $d=3$
is bounded by a constant independent of $h$: there exists a positive constant $\sigma$ independent of $h$ such that,
\begin{equation}
\label{eq:reg}
{
\ds \max_{\kappa \in \mathcal{T}_{h}} \frac{h_{ \kappa}}{\rho_{ \kappa}} \le \sigma.}
\end{equation}
\end{itemize}
As usual, $h$ denotes the maximal diameter
of all elements of $\mathcal{T}_{h}$.   To define the finite element functions,  let $r$ be a non negative integer.  For each {  $\kappa$ in $\mathcal{T}_{h}$, we denote by  $\mathbb{P}_r(\kappa)$  the space of restrictions to $\kappa$ of polynomials in $d$ variables and total degree at most  $r$, with a similar notation on the faces or edges of $\kappa$}. For every edge (when $d=2$) or face (when $d=3$) $e$ of the mesh  $\mathcal{T}_h$, we denote by $h_e$ the diameter of $e$.
%

\noindent {  We shall use the following inverse   inequality {  \cite{sayahdib}}:  for any dimension $d$,  there exists a constant $C_I$ such that for any polynomial function $v_h$ {  of degree $r$   on $\kappa$},
\begin{equation}
\label{eq:inversin}
 {  \|v_h\|_{L^3(\kappa)} \leq  C_I  h_\kappa^{-\frac{d}{6}}\|v_h\|_{L^2(\kappa)}.}
\end{equation}
The constant $C_I$ depends on the regularity parameter $\sigma$ of \eqref{eq:reg}, but for the sake of simplicity this is not indicated.}
\subsection{Discretization of the first variational problem}
In this section we follow the discretizations introduced in \cite{GW,LMJ,SJ}. In fact, the authors in \cite{GW} introduce a discrete variational formulation corresponding to \eqref{V1} based on the piecewise constant velocities and nonconforming piecewise $\P_1$ pressures (Crouzeix-Raviart element). They prove the existence and uniqueness of the exact and discrete solutions, and propose an alternating-directions algorithm inspired by the Peaceman-Rachford algorithm (see \cite{PR}), for solving the nonlinear system. In \cite{LMJ}, the authors carried out numerical tests of the methods studied in \cite{GW} and propose another mixed finite element spaces in which the approximation to the pressure is more smooth. Next, the work \cite{SJ} presented theoretical and numerical studies of the mixed finite element space, as proposed in \cite{LMJ}.

Let $X_{u,h} \subset X_u$ and $M_{u,h}\subset M_u$ the discrete spaces corresponding to the velocity and the pressure. We assume that they satisfy the following inf-sup condition:
\begin{equation}
\label{infsuph1}
\forall\, q_h\in M_{u,h},\; \sup_{\v_h\in X_{u,h}}\ds \frac{\ds \int_{\Omega} \nabla q_h \cdot \v_h \,  d\x\,}{\|\v_h\|_{X_u}}\geq \beta_u \|q_h\|_{M_{u,h}},
\end{equation}
where $\beta_u$ is a positive constant independent of $h$.\\

Problem \eqref{V1} can be discretized as following:
\begin{equation}\label{V1h}
\left\{
\begin{array}{ll}
\medskip
\forall \v_h \in X_{u,h}, \quad \ds \frac{\mu}{\rho} \int_\Omega K^{-1} \u_h \cdot \v_h \, d\x + \frac{\beta}{\rho} \int_\Omega |\u_h| \u_h \cdot \v_h \, d\x + \int_\Omega \nabla p_h \cdot \v_h \, d\x = \ds \int_\Omega \f \cdot \v_h \, d\x,\\
\forall q_h\in M_{u,h}, \quad \ds \int_\Omega \nabla q_h \cdot \u_h \, d\x = -\int_\Omega b q_h \, d\x + \int_\Gamma g_u q_h \, ds.
\end{array}
\right.
\end{equation}
In the following, we will consider for instance the discrete spaces introduced in \cite{SJ} given by:
\begin{equation}\label{defespaceprem}
\begin{split}
X_{u,h}=&\{  \v_h\in L^2(\bar{\Omega})^d; \, \forall\, \kappa \in \mathcal{T}_h, \; \v_h|_\kappa \in \P_0^d \}, \\
M_{u,h}=&\{  q_h\in C^0(\bar{\Omega}); \, \forall\, \kappa \in \mathcal{T}_h, \;q_{h}|_\kappa \in \P_1 \} \cap L^2_0(\Omega), \\
V_{u,h} =& \{ \v_h \in X_{u,h}; \forall q_h \in M_{u,h}, \ds \int_\Omega \nabla q_h \cdot \v_h \, d\x= 0  \}.
\end{split}
\end{equation}
It is shown in \cite{SJ} that there exists a unique $\u_{h,l} \in X_{u,h}$ such that
\begin{equation}\label{equat11}
\forall q_h \in M_{u,h}, \quad \ds \int_\Omega \nabla q_h \cdot \u_{h,l} \, d\x = \ds - \int_\Omega b\, q_h d\x + \int_\Gamma g \, q_h ds,
\end{equation}
and $\u_{h,l}$ verifies the following bound,
\begin{equation}\label{equat122222}
||\u_{h,l}||_{L^3(\Omega)^d/V_{u,h}} \le C_{l1} \big( ||b||_{L^{\frac{3d}{3+d}}(\Omega)} + ||g_u||_{L^{\frac{3(d-1)}{d}}(\Gamma)} \big).
\end{equation}
It is also shown in \cite{SJ} that Problem \eqref{V1h} admits a unique solution $(\u_h,p_h) \in X_{u,h} \times M_{u,h}$ satisfying exactly similar bounds as \eqref{bound1u}. Also the solutions $(\u,p)$ of \eqref{V1} and $(\u_h,p_h)$ of $\eqref{V1h}$ verify the following {\it a priori} error: \\
If $(\u,p)\in W^{1,4}(\Omega)^d \times W^{2,3/2}(\Omega)$, then there exists a constant $C$ independent of $h$ such that
\begin{equation}
\begin{array}{ll}
||\u - \u_h||_{L^2(\Omega)}  \le C h, \\
||\nabla (p - p_h) ||_{3/2,\Omega} \le C h.
\end{array}
\end{equation}

{\bf An iterative algorithm:} {  In order to approximate} the solution of the non-linear problem \eqref{V1h}, we introduce the following iterative algorithm: {  For} a given initial guess $\u_h^0 \in X_{u,h}$ and having $\u_h^i$ at each iteration $i$, we compute $(\u_h^{i+1}, p^{i+1}_h)$ solution of
\begin{equation}\label{V1hi}
\left\{
\begin{array}{ll}
\medskip
\forall \v_h \in X_{u,h}, \quad \ds \int_\Omega \alpha (\u_h^{i+1} - \u_h^i) \cdot \v_h \, d\x +  \frac{\mu}{\rho} \int_\Omega K^{-1} \u^{i+1}_h \cdot \v_h \, d\x + \frac{\beta}{\rho} \int_\Omega |\u^i_h| \u^{i+1}_h \cdot \v_h \, d\x \\
\medskip
\hspace{5cm} + \ds \int_\Omega \nabla p^{i+1}_h \cdot \v_h \, d\x = \ds \int_\Omega \f \cdot \v_h \, d\x,\\
\forall q_h\in M_{u,h}, \quad \ds \int_\Omega \nabla q_h \cdot \u^{i+1}_h \, d\x = -\int_\Omega b q_h \, d\x + \int_\Gamma g_u q_h \, ds,
\end{array}
\right.
\end{equation}
where $\alpha$ is {  a given} positive parameter.\\
In the following, we investigate the convergence of Scheme (\ref{V1hi}).
{\thm \label{boundu1} {  Problem $(\ref{V1hi})$ admits a unique solution $(\u_h^{i+1}, p_h^{i+1})$ in $X_{u,h} \times M_{u,h}$}. Furthermore, if the initial value $\u_h^0$ satisfies  the condition
\begin{equation}\label{cond1f}
|| \u_h^0 ||^2_{L^2(\Omega)} \le \ds L_1(\f,\u_{h,l}),
\end{equation}
where
\[
L_1(\f,\u_{h,l}) = \ds \frac{\mu K_m}{2\rho} \Big(  \ds  (\frac{3\rho}{2\mu K_m}  + \frac{1}{2}) ||\f||^2_{L^2(\Omega)} + (\frac{1}{2} + \alpha + \frac{2\mu K_M^2}{3\rho K_m} )  ||\u_{h,l} ||^2_{L^2(\Omega)}   +\frac{4\beta}{3\rho} || \u_{h,l} ||^3_{L^3(\Omega)}\Big),
\]
and if $\alpha$ satisfies the condition
$$
\ds \frac{\alpha}{4} > \ds  \frac{3\beta^2}{2\rho \mu K_m} {  C_I^4} h^{-\frac{2d}{3}}|| \u_{h,l} ||^2_{L^3(\Omega)}
+  \frac{4 \beta}{3 \rho}  {  C_I^6} h^{-d}  L_2(\f, L_1(\f,\u_{h,l})),
$$
where
$$  { L_2(\f, \eta) = \ds \frac{1}{\sqrt{\alpha}} \Big( \frac{3 \rho}{2 \mu K_m} \| \f \|^2_{L^2(\Omega)}
+ \frac{3 \mu K_M^2}{2 \rho K_m} \eta + \frac{3 \beta^2}{2 \mu \rho K_m}  C_I^6  h^{-d} \eta^2 \Big)^{1/2},} $$
then the solution of Problem $(\ref{V1hi})$ satisfies the estimates
\begin{equation}\label{equatt2}
|| \u_h^{i+1} ||^2_{L^2(\Omega)} \le \ds L_1(\f,\u_{h,l})
\end{equation}
and
\begin{equation}\label{equatt2a}
|| \u_h^{i+1} ||^3_{L^3(\Omega)} \le \ds \frac{3\rho}{\beta} \big( \frac{2\rho}{\mu K_m}   +\frac{\alpha}{2}  \big) L_1(\f,\u_{h,l}).
\end{equation}
}
\noindent \textbf {Proof.} {  To prove the existence and  uniqueness of the solution of Problem $(\ref{V1hi})$ which is a square linear system in finite dimension, it suffices to show the uniqueness. For a given $\u_h^i$, let $(\u_{h1}^{i+1},p_{h1}^{i+1})$ and $(\u_{h2}^{i+1},p_{h2}^{i+1})$ two different solutions of Problem $(\ref{V1hi})$ and let $\w_h= \u_{h1}^{i+1} - \u_{h2}^{i+1}$ and $\xi_h = p_{h1}^{i+1} - p_{h2}^{i+1}$, then $(\w_h,\xi_h)$ is the solution of the following problem:
\begin{equation*}
\left\{
\begin{array}{ll}
\medskip
\forall \v_h \in X_{u,h}, \quad \ds \int_\Omega \alpha \w_h \cdot \v_h \, d\x +  \frac{\mu}{\rho} \int_\Omega K^{-1} \w_h \cdot \v_h \, d\x + \frac{\beta}{\rho} \int_\Omega |\u^i_h| \w_h \cdot \v_h \, d\x + \ds \int_\Omega \nabla \xi_h \cdot \v_h \, d\x = 0,\\
\forall q_h\in M_{u,h}, \quad \ds \int_\Omega \nabla q_h \cdot \w_h \, d\x = 0.
\end{array}
\right.
\end{equation*}
By taking $(\v_h,q_h)=(\w_h,\xi_h)$ and by remarking that $\ds \frac{\beta}{\rho} \int_\Omega |\u^i_h| |\w_h|^2 \, d\x \ge 0$, we obtain by using the properties of $K^{-1}$ the following bound:
\[
(\alpha +  \ds \frac {K_m \mu}{\rho}) ||\w_h||^2_{L^2(\Omega)}  \le 0.
\]
Thus, we deduce that $\w_h = 0$. The inf-sup condition \eqref{infsuph1} deduces that $\xi_h=0$ and then, we get  the uniqueness of the solution of Problem $(\ref{V1hi})$.\\
}
%
\noindent To prove the bound \eqref{equatt2}, we need first to bound the error $\|   \u_h^{i+1} - \u_h^i \|_{L^2(\Omega)}$ with respect to the previous value $\u_h^{i}$. The second equation of Problem $(\ref{V1hi})$ allows us to deduce the relation
\begin{equation}\label{reluiuip1}
\forall q_h \in M_{u,h}, \quad \ds \int_\Omega \nabla q_h \cdot (\u^{i+1}_h - \u_h^i) \, d\x  = 0.
\end{equation}
Then, the  first equation of (\ref{V1hi}) with $\v_h = \u_h^{i+1} - \u_h^i$ gives
\[
\ds \alpha \| \u_{h}^{i+1} -\u_{h}^i\|^2_{L^2(\Omega)} + \frac{\mu}{\rho} \int_\Omega K^{-1} \u_{h}^{i+1} \cdot (\u_{h}^{i+1} -\u_{h}^i) \, d\x + \frac{\beta}{\rho} \int_\Omega |\u_{h}^i|  \u_{h}^{i+1} \cdot (\u_{h}^{i+1} -\u_{h}^i) \, d\x =  \ds \int_\Omega \f \cdot (\u_{h}^{i+1} - \u_h^i) \, d\x.
\]
By inserting $\pm \u_h^i$ in the second and the third terms of the last equation we get
\[
\begin{array}{ll}
\medskip
\ds \alpha \| \u_{h}^{i+1} -\u_{h}^i\|^2_{L^2(\Omega)} + \frac{\mu}{\rho} \int_\Omega K^{-1} |\u_{h}^{i+1} -\u_{h}^i |^2 \, d\x +
\frac{\beta}{\rho} \int_\Omega |\u_{h}^i| |\u_{h}^{i+1} -\u_{h}^i |^2 \, d\x \\
\hspace{2cm} = \ds \int_\Omega \f \cdot (\u_{h}^{i+1} - \u_h^i) \, d\x - \frac{\mu}{\rho} \int_\Omega K^{-1} \u_{h}^i \cdot (\u_{h}^{i+1} -\u_{h}^i) \, d\x
- \frac{\beta}{\rho} \int_\Omega |\u_{h}^i|  \u_{h}^i \cdot (\u_{h}^{i+1} -\u_{h}^i) \, d\x.
\end{array}
\]
By using the properties of $K$, the Cauchy-Schwartz inequality and Relation \eqref{eq:inversin},
we get, by remarking that the third term of the last equation is non-negative, the bound
\[
\begin{array}{ll}
\medskip
\ds \alpha \| \u_{h}^{i+1} -\u_{h}^i\|^2_{L^2(\Omega)} + \frac{\mu K_m}{\rho} \| \u_{h}^{i+1} -\u_{h}^i\|^2_{L^2(\Omega)} \le
\| \f \|_{L^2(\Omega)} \| \u_{h}^{i+1} -\u_{h}^i\|_{L^2(\Omega)} \\
\hspace{4cm} + \ds \frac{\mu K_M}{\rho}  \|\u_{h}^i \|_{L^2(\Omega)} \|\u_{h}^{i+1} -\u_{h}^i\|_{L^2(\Omega)} +
\frac{\beta}{\rho} {  C_I^3 }  h^{-\frac{d}{2}}   \|\u_{h}^i\|^2_{L^2(\Omega)}  \|\u_{h}^{i+1} -\u_{h}^i\|_{L^2(\Omega)}.
\end{array}
\]
We apply the relation $ab \le \ds \frac{1}{2\varepsilon} a^2 + \frac{1}{2} \varepsilon b^2$ with $\varepsilon=\ds \frac{\mu K_m}{3 \rho}$
for each term of the {  right hand side} of the previous inequality, and we get
\[
\begin{array}{ll}
\medskip
\ds \alpha \| \u_{h}^{i+1} -\u_{h}^i\|^2_{L^2(\Omega)} + \frac{\mu K_m}{2\rho} \| \u_{h}^{i+1} -\u_{h}^i\|^2_{L^2(\Omega)} \\
\hspace{3cm} \le
\ds \frac{3 \rho}{2\mu K_m} \| \f \|^2_{L^2(\Omega)} + \frac{3 \mu K_M^2}{2\rho K_m} \|\u_{h}^i \|^2_{L^2(\Omega)}
+ \frac{3 \beta^2}{2\mu \rho K_m} {  C_I^6 } h^{-d} \|\u_{h}^i\|^4_{L^2(\Omega)}
\end{array}
\]
and then the following bound
\begin{equation}\label{bound11}
\ds  \| \u_{h}^{i+1} -\u_{h}^i\|_{L^2(\Omega)} \le L_2(\f, \| \u_h^i \|^2_{L^2(\Omega)}),
\end{equation}
where $$  { L_2(\f, \eta) = \ds \frac{1}{\sqrt{\alpha}} \Big( \frac{3 \rho}{2 \mu K_m} \| \f \|^2_{L^2(\Omega)}
+ \frac{3 \mu K_M^2}{2 \rho K_m} \eta + \frac{3 \beta^2}{2 \mu \rho K_m}  C_I^6  h^{-d} \eta^2 \Big)^{1/2}.}$$
Then, we are in position to show the relation \eqref{equatt2}. Property \eqref{equat11} allows us to deduce that the term
$\u^{i+1}_{h,0} = \u_h^{i+1} - \u_{h,l}$ is in $X_{u,h}$ and verifies
\begin{equation}\label{equat111}
\forall q_h \in M_{u,h}, \quad \ds \int_\Omega \nabla q_h \cdot \u^{i+1}_{h,0} \, d\x  = 0.
\end{equation}
We consider the first equation of (\ref{V1hi}) with $\v_h = \u^{i+1}_{h,0} = \u_h^{i+1} - \u_{h,l}$ and we obtain:
\[
\begin{array}{ll}
\medskip
\ds \alpha \int_\Omega (\u_{h}^{i+1} -\u_{h}^i)  \cdot  \u_{h}^{i+1} \, d\x +  \frac{\mu}{\rho} \int_\Omega K^{-1} |\u_{h}^{i+1}|^2 \, d\x + \frac{\beta}{\rho} \int_\Omega |\u_{h}^{i+1}|^3  \, d\x =  \ds \int_\Omega \f \cdot (\u_{h}^{i+1} - \u_{hl}) \, d\x \\
\medskip
\hspace{1cm} \ds + \alpha \int_\Omega (\u_{h}^{i+1} -\u_{h}^i)  \cdot  \u_{h,l} \, d\x +  \frac{\mu}{\rho} \int_\Omega K^{-1} \u_{h}^{i+1} \cdot \u_{h,l}  \, d\x + \frac{\beta}{\rho} \int_\Omega (|\u_h^{i+1}| - |\u^i_h|) |\u_{h}^{i+1}|^2  \, d\x\\
\hspace{2cm} \ds + \frac{\beta}{\rho} \int_\Omega (|\u^i_h| - |\u^{i+1}_h|)\u_{h}^{i+1} \cdot \u_{h,l} \, d\x + \frac{\beta}{\rho} \int_\Omega |\u^{i+1}_h| \,\u_{h}^{i+1} \cdot \u_{h,l} \, d\x.
\end{array}
\]
By using the properties of $K$, the Cauchy-Schwartz inequality
and the relation $a^2 b \le \frac{1}{3} \big( \frac{1}{\delta^3} b^3 + 2 \delta^{3/2} a^3 \big)$ (for any positive real numbers $a$ and $b$),
we get:
\[
\begin{array}{ll}
\medskip
\ds \frac{\alpha}{2} ||\u_{h}^{i+1}||^2_{L^2(\Omega)} - \frac{\alpha}{2} ||\u_{h}^i||^2_{L^2(\Omega)} + \frac{\alpha}{2} ||\u_{h}^{i+1} -\u_{h}^i||^2_{L^2(\Omega)} +  \frac{\mu}{\rho} K_m ||\u_{h}^{i+1}||^2_{L^2(\Omega)} + \frac{\beta}{\rho}  || \u_{h}^{i+1} ||^3_{L^3(\Omega)}\\
\medskip
\hspace{1cm} \le \ds ||\f||_{L^2(\Omega)} ||\u_{h}^{i+1} ||_{L^2(\Omega)} + ||\f||_{L^2(\Omega)} ||\u_{h,l}||_{L^2(\Omega)} + \alpha ||\u_{h}^{i+1} -\u_{h}^i||_{L^2(\Omega)} ||\u_{h,l} ||_{L^2(\Omega)}\\
\medskip
\hspace{1.5cm} + \ds \frac{\mu}{\rho} K_M || \u_{h}^{i+1} ||_{L^2(\Omega)} || \u_{h,l} ||_{L^2(\Omega)} + \frac{\beta}{\rho} || \u^i_h - \u^{i+1}_h ||_{L^3(\Omega)} || \u_{h}^{i+1} ||_{L^3(\Omega)} || \u_{h,l} ||_{L^3(\Omega)}  \\
\hspace{2cm}  \ds  + \frac{\beta}{\rho} || \u^i_h - \u^{i+1}_h ||_{L^3(\Omega)} || \u_{h}^{i+1} ||^2_{L^3(\Omega)} + \frac{\beta}{\rho} ||\u^{i+1}_h ||^2_{L^3(\Omega)}  || \u_{h,l} ||_{L^3(\Omega)}.
\end{array}
\]
We deduce by using the relation \eqref{eq:inversin}, that for any positive numbers $\varepsilon_i,i=1 \dots,4$ and $\delta_j,j=1,2$, we have the following bound:
\[
\begin{array}{ll}
\medskip
\ds \frac{\alpha}{2} ||\u_{h}^{i+1}||^2_{L^2(\Omega)} - \frac{\alpha}{2} ||\u_{h}^i||^2_{L^2(\Omega)} + \frac{\alpha}{2} ||\u_{h}^{i+1} -\u_{h}^i||^2_{L^2(\Omega)} +  \frac{\mu}{\rho} K_m ||\u_{h}^{i+1}||^2_{L^2(\Omega)} + \frac{\beta}{\rho} || \u_{h}^{i+1} ||^3_{L^3(\Omega)}\\
\medskip
\hspace{1cm} \le \ds  \frac{1}{2 \varepsilon_1} ||\f||^2_{L^2(\Omega)} + \frac{1}{2} \varepsilon_1 ||\u_{h}^{i+1} ||^2_{L^2(\Omega)} + \frac{1}{2} ||\f||^2_{L^2(\Omega)} + \frac{1}{2} ||\u_{h,l} ||^2_{L^2(\Omega)}\\
\medskip
\hspace{1.5cm}  \ds + \frac{\alpha}{2 \varepsilon_2}  ||\u_{h}^{i+1} -\u_{h}^i||^2_{L^2(\Omega)} + \frac{\alpha}{2} \varepsilon_2 ||\u_{h,l}||^2_{L^2(\Omega)} + \frac{\mu^2}{2\rho^2 \varepsilon_3} K_M^2 || \u_{h,l} ||^2_{L^2(\Omega)} + \frac{1}{2} \varepsilon_3 || \u_{h}^{i+1} ||^2_{L^2(\Omega)}\\
\medskip
\hspace{1.5cm}  \ds + \frac{\beta^2}{2\rho^2 \varepsilon_4} {  C_I^4} h^{-\frac{2d}{3}}|| \u_{h,l} ||^2_{L^3(\Omega)} || \u^i_h - \u^{i+1}_h ||^2_{L^2(\Omega)} + \frac{1}{2} \varepsilon_4   || \u_{h}^{i+1} ||^2_{L^2(\Omega)} \\
\medskip
\hspace{1.5cm}  \ds +  \frac{\beta}{3\rho} \big(  (\frac{1}{\delta_1})^{3} {  C_I^6} h^{-d} ||\u^{i+1}_h -\u_{h}^i ||^3_{L^2(\Omega)}   +  2 \delta_1^{3/2} || \u_{h}^{i+1} ||^3_{L^3(\Omega)} \big)   \\

\hspace{1.5cm}  \ds + \frac{\beta}{3\rho} \big( (\frac{1}{\delta_2})^{3}  || \u_{h,l} ||^3_{L^3(\Omega)}   +  2 \delta_2^{3/2}||\u^{i+1}_h ||^3_{L^3(\Omega)} \big)
\end{array}
\]
We choose $\varepsilon_1=\varepsilon_3=\varepsilon_4=\ds\frac{\mu K_m}{3\rho}$, $\varepsilon_2=2$, $\delta_1=\delta_2=\ds \big(\frac{1}{2}\big)^{2/3}$ and we denote
$$
C_1 (\| \u_h^i \|_{L^2(\Omega)}) = \ds \frac{\alpha}{4} -  \frac{3\beta^2}{2\rho \mu K_m} {  C_I^4} h^{-\frac{2d}{3}}|| \u_{h,l} ||^2_{L^3(\Omega)}
-  \frac{4 \beta}{3 \rho}  {  C_I^6} h^{-d} L_2(\f, \| \u_h^i \|^2_{L^2(\Omega)}),
$$
which is not necessarily positive at this level.\\
By using the bound \eqref{bound11} , we get
$$
C_1 (\| \u_h^i \|_{L^2(\Omega)}) \le \ds \frac{\alpha}{4} -  \frac{3\beta^2}{2\rho \mu K_m} {  C_I^4} h^{-\frac{2d}{3}}|| \u_{h,l} ||^2_{L^3(\Omega)}
-  \frac{4 \beta}{3 \rho}  {  C_I^6}  h^{-d} \| \u_h^{i+1} - \u_h^i \|_{L^2(\Omega)}
$$
{  and} then we conclude the following bound
\begin{equation}\label{relata1}
\begin{array}{ll}
\medskip
\ds \frac{\alpha}{2} ||\u_{h}^{i+1}||^2_{L^2(\Omega)} - \frac{\alpha}{2} ||\u_{h}^i||^2_{L^2(\Omega)} + C_1(\| \u_h^i \|_{L^2(\Omega)}) ||\u_{h}^{i+1} -\u_{h}^i||^2_{L^2(\Omega)} +  \frac{\mu K_m}{2 \rho}||\u_{h}^{i+1}||^2_{L^2(\Omega)} + \frac{\beta}{3\rho} || \u_{h}^{i+1} ||^3_{L^3(\Omega)}\\
\medskip
\hspace{1cm} \le \ds  (\frac{3\rho}{2\mu K_m}  + \frac{1}{2}) ||\f||^2_{L^2(\Omega)} + (\frac{1}{2} + \alpha + \frac{2\mu K_M^2}{3\rho K_m} )  ||\u_{h,l} ||^2_{L^2(\Omega)}   +\frac{4\beta}{3\rho} || \u_{h,l} ||^3_{L^3(\Omega)}\\
\hspace{1cm} \le \ds \frac{2\rho }{\mu K_m} L_1(\f,\u_{h,l}).
\end{array}
\end{equation}
We now prove Estimate (\ref{equatt2}) by induction on $i$ under some conditions on $\alpha$. Starting with the relation \eqref{cond1f}, we suppose that we have
\begin{equation}\label{induction1}
|| \u_h^i ||^2_{L^2(\Omega)} \le \ds L_1(\f,\u_{h,l}).
\end{equation}
We are in one of the following two situations :
\begin{itemize}
\item We have $\ds  ||\u_h^{i+1}||_{L^2(\Omega)} \le
||\u_h^i||_{L^2(\Omega)}$. We obviously deduce the bound
\[
||\u_h^{i+1}||^2_{L^2(\Omega)} \le \ds L_1(\f,\u_{h,l})
\]
from the induction hypothesis.
\item We have $\ds  ||\u_h^{i+1}||_{L^2(\Omega)} \ge
||\u_h^i||_{L^2(\Omega)}$. By using the induction condition \eqref{induction1}, we chose
$$
\begin{array}{rcl}
\medskip
\ds \frac{\alpha}{4} &>& \ds  \frac{3\beta^2}{2\rho \mu K_m} {  C_I^4} h^{-\frac{2d}{3}}|| \u_{h,l} ||^2_{L^3(\Omega)}
+  \frac{4 \beta}{3 \rho}  {  C_I^6} h^{-d}  L_2(\f, L_1(\f,\u_{h,l}))\\
&>&  \ds  \frac{3\beta^2}{2\rho \mu K_m} {  C_I^4} h^{-\frac{2d}{3}}|| \u_{h,l} ||^2_{L^3(\Omega)}
+  \frac{4 \beta}{3 \rho}  {  C_I^6} h^{-d}  L_2(\f,||\u_h^i||^2_{L^2(\Omega)}  ),
\end{array}
$$
to get $C_1(||\u_h^i||_{L^2(\Omega)} ) >0$ and then to deduce that
\[
||\u_h^{i+1}||^2_{L^2(\Omega)} \le \ds L_1(\f,\u_{h,l}).
\]
\end{itemize}
whence we deduce the relation \eqref{equatt2}. The bound \eqref{equatt2a} is a simple consequence of Equation \eqref{relata1} and Relation \eqref{equatt2}. $\hfill\Box$\\

The next theorem show the convergence of the solution $(\u_h^i, p_h^i)$ of Problem $(\ref{V1hi})$ in $L^2(\Omega)^d \times L^2(\Omega)$ to the solution
$(\u_h, p_h)$ of Problem $(\ref{V1h})$.
{\thm \label{converg1} Assume that there exists ${  \beta_0 >0}$  such that, for every element {  $\kappa \in {\mathcal T}_h$, we have
\[ h_\kappa \ge \beta_0 h,\]}
(which means that the family of triangulations is uniformly regular). Under the assumptions of Theorem $\ref{boundu1}$ and if $\alpha$ satisfies also the condition
\begin{equation}\label{condalpha}
\ds \alpha > \ds C_1 h^{-d},
\end{equation}
where $$C_1=\ds \frac{\rho C^2}{K_m \mu}\qquad \mbox{and} \qquad C=\ds \frac{\beta}{\rho} {  C_I^3} (L_1(\f,\u_{h,l}))^{1/2},$$ then the sequence of solutions $(\u_h^i, p_h^i)$ of Problem $(\ref{V1hi})$ converges in $L^2(\Omega)^d \times L^2(\Omega)$ to the solution
$(\u_h, p_h)$ of Problem $(\ref{V1h})$. }

\noindent \textbf {Proof.} We take the difference between the
equations (\ref{V1hi}) and (\ref{V1h}) with $\v_h = \u^{i+1}_h -\u_h$ and we obtain the equation
\[
\begin{array}{rr}
\medskip
\ds \frac{\alpha}{2} || \u_h^{i+1} - \u_h ||^2_{L^2(\Omega)}  - \frac{\alpha}{2} || \u_h^{i} - \u_h ||^2_{L^2(\Omega)} +
\frac{\alpha}{2} || \u_h^{i+1} - \u^i_h ||^2_{L^2(\Omega)} + \frac{\mu}{\rho} \int_\Omega K^{-1} ( \u_h^{i+1} - \u_h)^2 \, d\x \\
+ \ds \frac{\beta}{\rho} (|\u_h^{i}| \u_h^{i+1} - |\u_h| \u_h , \u_h^{i+1} - \u_h) = 0.
\end{array}
\]
The last term in the previous equation, denoted by $T$, can be decomposed
as
\[
T = \ds \frac{\beta}{\rho} ((|\u_h^{i}| - |\u^{i+1}_h|)  \u_h^{i+1}, \u_h^{i+1} - \u_h) + \frac{\beta}{\rho} (|\u_h^{i+1}| \u_h^{i+1} -
|\u_h| \u_h , \u_h^{i+1} - \u_h).
\]
We denote by $T_1$ and $T_2$, respectively  the first and the
second terms in the right-hand side of the  last equation. Using
Lemma \ref{relatinf}, we have $T_2 \ge 0$. Then we derive by
using (\ref{eq:inversin}) and \eqref{KmM},
\[
\begin{array}{ll}
\medskip
\ds \frac{\alpha}{2} || \u_h^{i+1} - \u_h ||^2_{L^2(\Omega)}  - \frac{\alpha}{2} || \u_h^{i} - \u_h ||^2_{L^2(\Omega)} +
\frac{\alpha}{2}|| \u_h^{i+1} - \u^i_h ||^2_{L^2(\Omega)} + \frac{K_m \mu}{\rho}  \| \u_h^{i+1} - \u_h \|^2_{L^2(\Omega)} + T_2 \le  | T_1 |\\
\medskip
\hspace{1cm} \le \ds \frac{\beta}{\rho} \int_\Omega |\u_h^{i+1} - \u_h^{i}| \, |\u_h^{i+1}| \, |\u_h^{i+1} - \u_h| d{\textbf x}\\
\medskip
\hspace{1cm} \le \ds \frac{\beta}{\rho} ||\u_h^{i+1} - \u_h^{i} ||_{L^3(\Omega)} \; ||\u_h^{i+1}||_{L^3(\Omega)} \, ||\u_h^{i+1} - \u_h ||_{L^3(\Omega)}\\
\medskip
\hspace{1cm} \le \ds \frac{\beta}{\rho} {  C_I^3} h^{-\frac{d}{2}}||\u_h^{i+1} - \u_h^{i} ||_{L^2(\Omega)} \; ||\u_h^{i+1}||_{L^2(\Omega)} \, ||\u_h^{i+1} - \u_h ||_{L^2(\Omega)}\\
\medskip
\hspace{1cm} \le \ds \frac{\beta}{\rho} {  C_I^3} h^{-\frac{d}{2}} (L_1(\f,\u_{h,l})^{1/2}  ||\u_h^{i+1} - \u_h^{i} ||_{L^2(\Omega)} \, ||\u_h^{i+1} - \u_h ||_{L^2(\Omega)}.
\end{array}
\]
We denote by $C= \ds \frac{\beta}{\rho} {  C_I^3} (L_1(\f,\u_{h,l})^{1/2}  $ and we use the inequality $ab \le \ds
\frac{1}{2\varepsilon} a^2 + \frac{\varepsilon}{2} b^2$ (with
$\varepsilon = \ds \frac{K_m \mu}{\rho}$) to obtain the
following bound
\[
\begin{array}{ll}
\medskip
\ds \frac{\alpha}{2} || \u_h^{i+1} - \u_h ||^2_{L^2(\Omega)}  - \frac{\alpha}{2} || \u_h^{i} - \u_h ||^2_{L^2(\Omega)} +
\frac{\alpha}{2} || \u_h^{i+1} - \u^i_h ||^2_{L^2(\Omega)} + \frac{K_m \mu}{2\rho } {  || \u_h^{i+1} - \u_h ||^2_{L^2(\Omega)}} \\
\hspace{2cm} \le \ds \frac{\rho C^2}{2K_m \mu} h^{-d} || \u_h^{i+1} - \u_h^i ||^2_{L^2(\Omega)}.
\end{array}
\]
We {  choose} $\ds \frac{\alpha}{2} > \ds \frac{\rho C^2}{2K_m \mu} h^{-d}$, denote by $C_1 =
\ds \frac{1}{2}(\alpha - \ds \frac{\rho C^2}{K_m \mu} h^{-d})$  and obtain
\begin{equation}\label{relat22}
\ds \frac{\alpha}{2} || \u_h^{i+1} - \u_h ||^2_{L^2(\Omega)}  - \frac{\alpha}{2} || \u_h^{i} - \u_h ||^2_{L^2(\Omega)} + C_1 ||
\u_h^{i+1} - \u^i_h ||^2_{L^2(\Omega)} + \frac{\mu K_m}{2\rho} {  || \u_h^{i+1} - \u_h ||^2_{L^2(\Omega)} } \le 0.
\end{equation}
We deduce that, for all $i\ge 1$, we have (if $|| \u_h^{i} - \u_h ||_{L^2(\Omega)} \neq 0$)

\[
|| \u_h^{i+1} - \u_h ||_{L^2(\Omega)}  < || \u_h^{i} - \u_h
||_{L^2(\Omega)},
\]
and we deduce the convergence of the sequence $(\u_h^{i+1} - \u_h)$ in $L^2(\Omega)^d$ and then the convergence of the sequence $\u_h^{i}$ in $L^2(\Omega)^d$. By taking the limit of (\ref{relat22}) and remarking that $T_2 \ge 0$, we get
\[
\underset{{\small i \rightarrow +\infty} }{\lim} \Big( {  || \u_h^{i+1} - \u_h ||^2_{L^2(\Omega)}}  \Big) \le 0.
\]
We deduce then that ${  || \u_h^{i+1} - \u_h ||_{L^2(\Omega)} }$ converges to $0$ and $\u_h^{i+1}$ converges to $\u_h$ in $L^2(\Omega)^d$.\\
For the convergence of the pressure, We take the difference between the equations (\ref{V1hi}) and (\ref{V1h}) and we obtain for all $\v_h \in X_{u,h}$ the equation
\[
\begin{array}{rr}
\medskip
\ds \int_\Omega \nabla (p_h^{i+1} - p_h) \v_h \, d\x = \ds -\alpha \int_\Omega (\u_h^{i+1} - \u^i_h)\v_h \, d\x  + \frac{\mu}{\rho} \int_\Omega K^{-1} ( \u_h - \u_h^{i+1}) \v_h \, d\x \\
\hspace{2cm} \ds  + \frac{\beta}{\rho} ( (|\u_h| - |\u_h^i |) \u_h, \v_h) + \frac{\beta}{\rho} ( |\u^i_h| (\u_h - \u_h^{i+1}) , \v_h).
\end{array}
\]
We get by using the inverse inequality \eqref{eq:inversin} the following:
\[
\begin{array}{ll}
\medskip
\ds \frac{\Big| \int_\Omega \nabla (p_h^{i+1} - p_h) \v_h \, d\x \Big| }{||\v_h ||_{L^3(\Omega)}} \le  \ds (\alpha || \u_h^i - \u_h^{i+1}||_{L^2(\Omega)} + \frac{\mu K_m}{\rho}  || \u_h - \u_h^{i+1}||_{L^2(\Omega)}) \frac{||\v_h||_{L^2(\Omega)}}{||\v_h||_{L^3(\Omega)}} \\
 \hspace{2cm} \ds  { + \frac{\beta}{\rho}  C_I h^{-\frac{d}{6}}  ||\u_h - \u_h^i||_{L^2(\Omega)} (||\u_h ||_{L^3(\Omega)}  + ||\u_h^i||_{L^3(\Omega)} ).}
\end{array}
\]
Owning the inf-sup condition \eqref{infsuph1}, we deduce the relation
\[
\begin{array}{ll}
\medskip
\ds || \nabla (p_h^{i+1} - p_h) ||_{L^{3/2}(\Omega)} \le  \ds \frac{1}{\beta_u}\Big( \alpha {  |\Omega|^{1/6}} || \u_h^i - \u_h^{i+1}||_{L^2(\Omega)} + {  |\Omega|^{1/6}} \frac{\mu K_m}{\rho}  || \u_h - \u_h^{i+1}||_{L^2(\Omega)}\\
 \hspace{2cm} \ds {  + \frac{\beta}{\rho}  C_I h^{-\frac{d}{6}}  ||\u_h - \u_h^i||_{L^2(\Omega)} (||\u_h ||_{L^3(\Omega)}  + ||\u_h^i||_{L^3(\Omega)} )\Big).}
\end{array}
\]
{  Thus for a given mesh (given $h$), the strong convergence of $\u_h^i$ to $\u_h$ in $L^2(\Omega)^d$ deduces the strong convergence of $\nabla p_h^i$ to $\nabla p_h$ in $L^{\frac{3}{2}}(\Omega)$. Furthermore, the fact that $p_h^i$ and $p_h$ are in the discrete space of $\P_1$ finite elements  $M_{u,h}\subset L^2_0(\Omega)$ defined in \eqref{defespaceprem} allows us to deduce the strong convergence of $p_h^i$ to $p_h$ in $L^2(\Omega)$ .}
{\rmq The condition $(\ref{cond1f})$  supposes that the initial values of the algorithms are small related to the data. We can for example take $\u_h^0=\0$.}

{\rmq \label{remtheoralpha}
{  As I mentioned in the introduction, the works [5,7,12] (except the alternating-directions algorithm in [5]) present iterative schemes without studying the corresponding convergence. In this work, I introduce and study the convergence of the numerical scheme (3.9) which can be applied to all the discrete variational formulations presented in the above references. Theorems \ref{converg1} and \ref{boundu1} require that $\alpha > C h^{-d}$ to get the convergence of the numerical scheme (\ref{V1hi}), where the constant $C$ is not easy to compute practically, especially for the numerical investigations. In fact, This result of convergence states that for a given mesh (for a given $h$), the iterative solution $(u_h^i, p_h^i)$ converges to $(u_h,p_h)$ when $i \rightarrow +\infty$ for $\alpha > C h^{-d}$. Of course in this case, we can see clearly that this convergence depends on $h$. If we are interested to study this convergence when $h$ is very small or $h$ goes to $0$, we need to take $\alpha$ depending on $h$ in Scheme (\ref{V1hi}) such that it verifies the condition $\alpha > C h^{-d}$ introduced in Theorems \ref{converg1} and \ref{boundu1}.\\
Otherwise, Theorems \ref{converg1} and \ref{boundu1} show a necessary condition ($\alpha > C h^{-d}$) for the convergence, but perhaps it is not a sufficient one for the convergence. This is the best we can do till now and we are trying to get a more suitable convergence condition which will be less restrictive to this one.\\
Furthermore, we refer to Remark \ref{remrkconvnum} in Section \ref{sec4.1} for the discussion of the practical choice of $\alpha$ for the numerical investigations.
}

\subsection{Discretization of the second variational problem}
In this section we follow the discretization introduced in \cite{HH}. In fact, the authors in \cite{HH} introduce a discrete variational formulation corresponding to \eqref{V2} based on the mixed elements such as the Raviart-Thomas mixed element \cite{RaTho1} and Brezzi-Douglas-Marini mixed element \cite{BDM}. They also introduce an iterative scheme without studying the corresponding convergence, and show numerical investigations testifying the convergence of the finite element approximation. \\

In this section, we consider the discrete variational formulation, introduce a new corresponding numerical algorithm and show the convergence of the corresponding iterative solution. \\
Let $X_{p,h} \subset X_p$ and $M_{p,h}\subset M_p$ the discrete spaces corresponding to the velocity and the pressure. We assume that they satisfy the following inf-sup condition:
\begin{equation}
\label{infsuph2}
\forall\, q_h\in M_{p,h},\; \sup_{\v_h\in X_{p,h}}\ds \frac{\ds \int_{\Omega} q_h \, \div(\v_h) \,  d\x\,}{\|\v_h\|_{X_p}}\geq \beta_p \|q_h\|_{M_{p,h}},
\end{equation}
where $\beta_p$ is a positive constant independent of $h$.\\

Problem \eqref{V2} can be discretized as following:
\begin{equation}\label{V2h}
\left\{
\begin{array}{ll}
\medskip
\forall \v_h \in X_{p,h}, \quad \ds \frac{\mu}{\rho} \int_\Omega K^{-1} \u_h \cdot \v_h \, d\x + \frac{\beta}{\rho} \int_\Omega |\u_h| \u_h \cdot \v_h \, d\x - \int_\Omega  p_h \, \div(\v_h) \, d\x =
\ds \int_\Omega \f \cdot \v_h \, d\x,\\
\forall q_h\in M_{p,h}, \quad \ds \int_\Omega  q_h \, \div(\u_h) \, d\x = \ds \int_\Omega b q_h \, d\x.
\end{array}
\right.
\end{equation}
In the following, we will consider for instance the {  Raviart-Thomas $RT0$ mixed element \cite{RaTho1}} given by:
\begin{equation}
\begin{split}
X_{p,h}=& {  \{\v_h\in X_p; \, \, \v_{h}(\x)|_\kappa=a_\kappa \x +\b_\kappa, a_\kappa \in \R,\b_\kappa \in \R^d,\,\forall\, \kappa \in \mathcal{T}_h\},}\\
M_{p,h}=&\{q_h\in L^2(\Omega);\;\forall\, \kappa \in \mathcal{T}_h, \;q_{h}|_\kappa \, \mbox{is constant}\}.\\
\end{split}
\end{equation}
It is also shown in \cite{HH} that there exists a unique $\u_{h,p} \in X_{p,h}$ such that
\begin{equation}\label{equat12}
\forall q_h \in M_{p,h}, \quad \ds \int_\Omega q_h \, \div(\u_{h,p}) \, d\x = \ds  \int_\Omega b\, q_h d\x,
\end{equation}
and $\u_{h,p}$ verifies the following bound,
\begin{equation}\label{equat22}
||\u_{h,p}||_{X_{p}} \le C_{l2} ||b||_{L^2(\Omega)}.
\end{equation}
It is shown in \cite{HH} that Problem \eqref{V2h} admits a unique solution $(\u_h,p_h) \in X_{p,h} \times M_{p,h}$ satisfying exactly similar bounds as \eqref{bound2u}. Also the solutions $(\u,p)$ of \eqref{V2} and $(\u_h,p_h)$ of $\eqref{V2h}$ verify the following {\it a priori} error: \\
If $(\u,p)\in W^{1,3}(\Omega)^d \times W^{1,3/2}(\Omega)$, then there exists a constant $C$ independent of $h$ such that
\begin{equation}
\begin{array}{ll}
||\u - \u_h||_{L^2(\Omega)} + ||\u - \u_h||^3_{3,\Omega} \le C h^2,\\
||p - p_h ||_{L^2(\Omega)} \le C h.
\end{array}
\end{equation}

{\bf An iterative algorithm:} {  In order to approximate} the solution of the non-linear problem \eqref{V2h}, we introduce the following iterative algorithm: for a given initial guess $\u_h^0 \in X_{p,h}$ and having $\u_h^i$ at each iteration $i$, we compute $(\u_h^{i+1}, p^{i+1}_h)$ solution of
\begin{equation}\label{V2hi}
\left\{
\begin{array}{ll}
\medskip
\forall \v_h \in X_{p,h}, \quad \ds \alpha \int_\Omega (\u_h^{i+1} - \u_h^i) \v_h +  \frac{\mu}{\rho} \int_\Omega K^{-1} \u^{i+1}_h \cdot \v_h \, d\x + \frac{\beta}{\rho} \int_\Omega |\u^i_h| \u^{i+1}_h \cdot \v_h \, d\x \\
\medskip
\hspace{5cm} - \ds \int_\Omega  p^{i+1}_h \, \div(\v_h) \, d\x = \ds \int_\Omega \f \cdot \v_h \, d\x,\\
\forall q_h\in M_{p,h}, \quad \ds \int_\Omega  q_h \, \div(\u^{i+1}_h) \, d\x = \int_\Omega b q_h \, d\x,
\end{array}
\right.
\end{equation}
where $\alpha$ is {  a given} positive parameter.\\
In the following, we investigate the convergence of the scheme $(\ref{V2hi})$.
{\thm\label{boundu2} {  Problem $(\ref{V2hi})$ admits a unique solution $(\u_h^{i+1}, p_h^{i+1})$ in $X_{p,h} \times M_{p,h}$}. Furthermore, if the initial value $\u_h^0$ satisfies  the condition
\begin{equation}\label{cond1f1}
|| \u_h^0 ||^2_{L^2(\Omega)} \le \ds L_1(\f,\u_{h,p}),
\end{equation}
where
\[
L_1(\f,\u_{h,p}) = \ds \frac{\mu K_m}{2\rho} \Big(  \ds  (\frac{3\rho}{2\mu K_m}  + \frac{1}{2}) ||\f||^2_{L^2(\Omega)}
+ (\frac{1}{2} + \alpha + \frac{2\mu K_M^2}{3\rho K_m} )  ||\u_{h,p} ||^2_{L^2(\Omega)}   +\frac{4\beta}{3\rho} || \u_{h,p} ||^3_{L^3(\Omega)}\Big),
\]
and if $\alpha$ satisfies the condition
$$
\ds \frac{\alpha}{4} > \ds  \frac{3\beta^2}{2\rho \mu K_m} {  C_I^4} h^{-\frac{2d}{3}}|| \u_{h,p} ||^2_{L^3(\Omega)}
+  \frac{4 \beta}{3 \rho}  {  C_I^6} h^{-d}  L_2(\f, L_1(\f,\u_{h,p})),
$$
where
$$  { L_2(\f, \eta) = \ds \frac{1}{\sqrt{\alpha}} \Big( \frac{3 \rho}{2 \mu K_m} \| \f \|^2_{L^2(\Omega)}
+ \frac{3 \mu K_M^2}{2 \rho K_m} \eta + \frac{3 \beta^2}{2 \mu \rho K_m}  C_I^6  h^{-d} \eta^2 \Big)^{1/2},} $$
then the solution of Problem $(\ref{V2hi})$ satisfies the estimates
\begin{equation}\label{equatt22}
|| \u_h^{i+1} ||^2_{L^2(\Omega)} \le \ds L_1(\f,\u_{h,p})
\end{equation}
and
\begin{equation}\label{equatt222a}
|| \u_h^{i+1} ||^3_{L^3(\Omega)} \le \ds \frac{3\rho}{\beta} \big( \frac{2\rho}{\mu K_m}   +\frac{\alpha}{2}  \big) L_1(\f,\u_{h,p}).
\end{equation}
}
\noindent \textbf {Proof.} {  Problem $(\ref{V2hi})$ is a square linear system in finite dimension. Then to prove the existence and uniqueness of the corresponding solution, it suffices to prove the uniqueness. For a given $\u_h^i$, let $(\u_{h1}^{i+1},p_{h1}^{i+1})$ and $(\u_{h2}^{i+1},p_{h2}^{i+1})$ two different solutions of Problem $(\ref{V2hi})$ and $\w_h= \u_{h1}^{i+1} - \u_{h2}^{i+1}$ and $\xi_h = p_{h1}^{i+1} - p_{h2}^{i+1}$, then $(\w_h,\xi_h)$ is the solution of the following problem:
\begin{equation*}
\left\{
\begin{array}{ll}
\medskip
\forall \v_h \in X_{p,h}, \quad \ds \alpha \int_\Omega \w_h \v_h +  \frac{\mu}{\rho} \int_\Omega K^{-1} \w_h \cdot \v_h \, d\x + \frac{\beta}{\rho} \int_\Omega |\u^i_h| \w_h \cdot \v_h \, d\x  - \int_\Omega \xi_h \div(\v_h) \, d\x = 0,\\
\forall q_h\in M_{p,h}, \quad \ds \int_\Omega  q_h \, \div(\w_h) \, d\x = 0,
\end{array}
\right.
\end{equation*}
By taking $(\v_h,q_h)=(\w_h,\xi_h)$ and by remarking that $\ds \frac{\beta}{\rho} \int_\Omega |\u^i_h| |\w_h|^2 \, d\x \ge 0$, we obtain by using the properties of $K^{-1}$ the following bound:
\[
(\alpha +  \ds \frac {K_m \mu}{\rho}) ||\w_h||^2_{L^2(\Omega)}  \le 0.
\]
Thus, we deduce that $\w_h = 0$. The inf-sup condition \eqref{infsuph2} deduces that $\xi_h=0$ and then, we get  the uniqueness of the solution of Problem $(\ref{V2hi})$.\\

\noindent Let us now prove the bound \eqref{equatt22}. We need first to bound the error $\|   \u_h^{i+1} - \u_h^i \|_{L^2(\Omega)}$ with respect to the previous value $\u_h^{i}$. The second equation of Problem $(\ref{V2hi})$ allows us to deduce the relation
\begin{equation*}
\forall q_h \in M_{u,h}, \quad \ds \int_\Omega  q_h \, \div(\u^{i+1}_h - \u_h^i) \, d\x  = 0.
\end{equation*}
Then, the  first equation of (\ref{V2hi}) with $\v_h = \u_h^{i+1} - \u_h^i$ gives
\[
\ds \alpha \| \u_{h}^{i+1} -\u_{h}^i\|^2_{L^2(\Omega)} + \frac{\mu}{\rho} \int_\Omega K^{-1} \u_{h}^{i+1} \cdot (\u_{h}^{i+1} -\u_{h}^i) \, d\x + \frac{\beta}{\rho} \int_\Omega |\u_{h}^i|  \u_{h}^{i+1} \cdot (\u_{h}^{i+1} -\u_{h}^i) \, d\x =  \ds \int_\Omega \f \cdot (\u_{h}^{i+1} - \u_h^i) \, d\x.
\]
This last equation can be treated exactly as its analogue one in the proof of Theorem \ref{boundu1} to get the following bound:
\begin{equation}\label{bound22}
\ds  \| \u_{h}^{i+1} -\u_{h}^i\|_{L^2(\Omega)} \le L_2(\f, \| \u_h^i \|^2_{L^2(\Omega)}),
\end{equation}
where $$  { L_2(\f, \eta) = \ds \frac{1}{\sqrt{\alpha}} \Big( \frac{3 \rho}{2 \mu K_m} \| \f \|^2_{L^2(\Omega)}
+ \frac{3 \mu K_M^2}{2 \rho K_m} \eta + \frac{3 \beta^2}{2 \mu \rho K_m}  C_I^6  h^{-d} \eta^2 \Big)^{1/2}.}$$
Now, Relation \eqref{bound22} allows us to show \eqref{equatt22}. In fact, Property \eqref{equat12} allows us to deduce that the term $\u^{i+1}_{h,0} = \u_h^{i+1} - \u_{h,p}$ is in $X_{p,h}$ and verifies
\begin{equation*}
\forall q_h \in M_{p,h}, \quad \ds \int_\Omega  q_h \, \div(\u^{i+1}_{h,0}) \, d\x  = 0.
\end{equation*}
We consider the first equation of (\ref{V2hi}) with $\v_h = \u^{i+1}_{h,0} = \u_h^{i+1} - \u_{h,p}$ and we obtain:
\[
\begin{array}{ll}
\medskip
\ds \alpha \int_\Omega (\u_{h}^{i+1} -\u_{h}^i)  \cdot  \u_{h}^{i+1} \, d\x +  \frac{\mu}{\rho} \int_\Omega K^{-1} |\u_{h}^{i+1}|^2 \, d\x + \frac{\beta}{\rho} \int_\Omega |\u_{h}^{i+1}|^3  \, d\x =  \ds \int_\Omega \f \cdot (\u_{h}^{i+1} - \u_{h,p}) \, d\x \\
\medskip
\hspace{1cm} \ds + \alpha \int_\Omega (\u_{h}^{i+1} -\u_{h}^i)  \cdot  \u_{h,p} \, d\x +  \frac{\mu}{\rho} \int_\Omega K^{-1} \u_{h}^{i+1} \cdot \u_{h,p}  \, d\x + \frac{\beta}{\rho} \int_\Omega (|\u_h^{i+1}| - |\u^i_h|) |\u_{h}^{i+1}|^2  \, d\x\\
\hspace{2cm} \ds + \frac{\beta}{\rho} \int_\Omega (|\u^i_h| - |\u^{i+1}_h|)\u_{h}^{i+1} \cdot \u_{h,p} \, d\x + \frac{\beta}{\rho} \int_\Omega |\u^{i+1}_h| \,\u_{h}^{i+1} \cdot \u_{h,p} \, d\x.
\end{array}
\]
From now on, the steps are identical to those used in the proof of Theorem \ref{boundu1} and it is useless to rewrite them. Then we get the bounds \eqref{equatt22} and \eqref{equatt222a}. }\\

The next theorem treats the convergence of the scheme \eqref{V2hi}.
{\thm Assume that there exists ${  \beta_0 >0}$  such that, for every element {  $\kappa\in {\mathcal T}_h$, we have
\[ h_\kappa \ge \beta_0 h,\]}
(which means that the family of triangulations is uniformly regular). Under the assumptions of Theorem $\ref{boundu2}$ and if $\alpha$ satisfies also the condition
\begin{equation}\label{condalpha2}
\ds \alpha > \ds C_1 h^{-d},
\end{equation}
where $$C_1=\ds \frac{\rho C^2}{K_m \mu}\qquad \mbox{and} \qquad C=\ds \frac{\beta}{\rho} {  C_I^2} L_1(\f,\u_{h,p})^{1/2},$$ then the sequence of solutions $(\u_h^i,p_h^i)$ of Problem $(\ref{V2hi})$ converges in $L^2(\Omega)^d \times L^2(\Omega)$ to the solution $(\u_h, p_h)$ of Problem $(\ref{V2h})$. }

\noindent \textbf {Proof.} We take the difference between the equations (\ref{V2h}) and (\ref{V2hi}) with $\v_h = \u^{i+1}_h -
\u_h$ and {  we follow similar steps} of the proof of Theorem \ref{converg1} to get the {  the convergence of $\u_h^i$ to $\u_h$ in $L^2(\Omega)^d$. In fact, it is useless to rewrite the corresponding details as they are strictly the same. \\
Now, we prove the convergence of the iterative pressure. We take the difference between the equations (\ref{V2h}) and (\ref{V2hi}) and we obtain for all $\v_h \in X_{p,h}$ the equation
\[
\begin{array}{rr}
\medskip
\ds \int_\Omega  (p_h^{i+1} - p_h) \div(\v_h) \, d\x = \ds \alpha \int_\Omega (\u_h^{i+1} - \u^i_h)\v_h \, d\x  - \frac{\mu}{\rho} \int_\Omega K^{-1} ( \u_h - \u_h^{i+1}) \v_h \, d\x \\
\hspace{2cm} \ds  - \frac{\beta}{\rho} ( (|\u_h| - |\u_h^i |) \u_h, \v_h) - \frac{\beta}{\rho} ( |\u^i_h| (\u_h - \u_h^{i+1}) , \v_h).
\end{array}
\]
We get by using the inverse inequality \eqref{eq:inversin} the following:
\[
\begin{array}{ll}
\medskip
\ds \frac{\Big| \int_\Omega  (p_h^{i+1} - p_h) \div(\v_h) \, d\x \Big| }{||\v_h ||_{X_p}} \le  \ds \big(\alpha || \u_h^i - \u_h^{i+1}||_{L^2(\Omega)} + \frac{\mu K_m}{\rho}  || \u_h - \u_h^{i+1}||_{L^2(\Omega)} \big) \frac{||\v_h||_{L^2(\Omega)}}{||\v_h||_{X_p}} \\
 \hspace{2cm} \ds  + \frac{\beta}{\rho} {  C_I} h^{-\frac{d}{6}}  ||\u_h - \u_h^i||_{L^2(\Omega)} (||\u_h ||_{L^3(\Omega)}  + ||\u_h^i||_{L^3(\Omega)} ) \frac{||\v_h||_{L^3(\Omega)}}{||\v_h||_{X_p}}.
\end{array}
\]
As we have $||\v_h||_{L^3(\Omega)} \le ||\v_h||_{X_p}$ and $||\v_h||_{L^2(\Omega)} \le |\Omega|^{1/6} ||\v_h||_{L^3(\Omega)}$, we deduce by using the inf-sup condition \eqref{infsuph2}, the following relation
\[
\begin{array}{ll}
\medskip
\ds || p_h^{i+1} - p_h ||_{L^2(\Omega)} \le  \ds \frac{1}{\beta_p}\Big( \alpha |\Omega|^{1/6} || \u_h^i - \u_h^{i+1}||_{L^2(\Omega)} + |\Omega|^{1/6}\frac{\mu K_m}{\rho}  || \u_h - \u_h^{i+1}||_{L^2(\Omega)}\\
 \hspace{2cm} \ds {  + \frac{\beta}{\rho}  C_I h^{-\frac{d}{6}}  ||\u_h - \u_h^i||_{L^2(\Omega)^d} (||\u_h ||_{L^3(\Omega)^d}  + ||\u_h^i||_{L^3(\Omega)^d} )\Big).}
\end{array}
\]
Thus, the strong convergence of $\u_h^i$ to $\u_h$ in $L^2(\Omega)^d$ deduces the strong convergence of $p_h^i$ to $p_h$ in $L^2(\Omega)$.
}
\section{Numerical results}

\noindent In this section, we present  numerical experiments for
our nonlinear problem. These simulations have been performed using
the code FreeFem++ due to F. Hecht and O. Pironneau, see
\cite{Freefem}.\\

We consider the domain $\Omega=]0,1[^2 \subset \R^2$, each edge is divided
into $N$ equal segments so that $\Omega$ is divided into $N^2$
equal squares and finally into $2N^2$ equal triangles. For simplicity, we take $\mu = \rho =1$ and $K=I$.\\

We will show in this section numerical investigations corresponding to Problems \eqref{V1hi} and \eqref{V2hi}.
In both cases, we use for the convergence the stopping criterion ${  Err_L}\leq \varepsilon$ where $\varepsilon$
is a given tolerance considered in this work equal to $10^{-5}$ and ${  Err_L}$ is defined by
$$Err_L=\ds \sqrt{\frac{||\u_h^{i+1}-\u_h^{i}||^2_{L^2(\Omega)} + ||p_h^{i+1}-p_h^{i}||^2_{L^2(\Omega)}}{||\u_h^{i+1}||^2_{L^2(\Omega)} + ||p_h^{i+1}||^2_{L^2(\Omega)}}}.$$
In both Algorithms \eqref{V1hi} and \eqref{V2hi}, the initial guess $\u^0_h$ {  is considered in one of this two situations:
\begin{enumerate}
\item $\u_h^0 = \0$.
\item $\u_h^0$ is calculated by using the Darcy's problem which corresponds to $\beta=\alpha=0$.
\end{enumerate}
We will see later that the second case where $\u_h^0$ is the solution of Darcy's problem improved the convergence of the algorithms.} \\
We consider also the error
\[
Err=\ds\sqrt{\frac{||\u_h^{i}-\u||^2_{L^2(\Omega)} + ||p_h^{i}-p||^2_{L^2(\Omega)}}{||\u||^2_{L^2(\Omega)} + ||p||^2_{L^2(\Omega)}}}
\]
which describes the convergence of the algorithms (\ref{V1hi}) and
(\ref{V2hi}).
\subsection{First discrete scheme \eqref{V1hi}}\label{sec4.1}
In this section we consider the scheme \eqref{V1hi} and show numerical corresponding tests. {  In fact, to compute the solution of the iterative problem \eqref{V1hi}, we use the penalty method (see \cite{Freefem}) which consists to solve the following problem:
\begin{equation*}
\left\{
\begin{array}{ll}
\medskip
\forall \v_h \in X_{u,h}, \quad \ds \int_\Omega \alpha (\u_h^{i+1} - \u_h^i) \cdot \v_h \, d\x +  \frac{\mu}{\rho} \int_\Omega K^{-1} \u^{i+1}_h \cdot \v_h \, d\x + \frac{\beta}{\rho} \int_\Omega |\u^i_h| \u^{i+1}_h \cdot \v_h \, d\x \\
\medskip
\hspace{5cm} + \ds \int_\Omega \nabla p^{i+1}_h \cdot \v_h \, d\x = \ds \int_\Omega \f \cdot \v_h \, d\x,\\
\forall q_h\in M_{u,h}, \quad \ds \int_\Omega \nabla q_h \cdot \u^{i+1}_h \, d\x - \varepsilon \int_\Omega p_h^{i+1} \, q_h \, d\x = -\int_\Omega b q_h \, d\x + \int_\Gamma g_u q_h \, ds,
\end{array}
\right.
\end{equation*}
where $\varepsilon$ is very small penalty coefficient. Here we choose $\varepsilon=10^{-8}$.\\
}

We consider the two following examples:
\begin{enumerate}
\item First example:
\begin{equation}\label{Ex1FA}
\left\{
\begin{array}{lcl}
\medskip
p(x,y) = \cos (\pi x) \cos (\pi y), &\qquad& \u(x, y) = \gamma (\exp(x) \, \sin (\pi y),  \ds \frac{1}{\pi} \exp(x) \, \cos (\pi y))^T ,\\
b=0, &\qquad&  \f =	  \u +  \beta |\u| \u+ \nabla p,
\end{array}
\right.
\end{equation}
and $\u\cdot\n=\ds -\frac{\gamma}{\pi} e^x \; \; \mbox{for} \;\; y=0,1;  \qquad \u\cdot\n= \gamma \sin(\pi y) \; \; \mbox{for} \;\; x=0;  \qquad \u\cdot\n= \gamma e \, \sin(\pi y) \; \; \mbox{for} \;\; x=1.$
\item Second example:
\begin{equation}\label{Ex2FA}
\left\{
\begin{array}{lcl}
\medskip
p(x, y) = x y^2-yx^2, &\qquad& \u(x, y) =  \gamma (x \exp(\pi y), y \exp(\pi x))^T ,\\
b=\gamma (\exp(\pi x)+ \exp(\pi y)), &\qquad&  \f =	  \u +  \beta |\u| \u+ \nabla p,
\end{array}
\right.
\end{equation}
and $\u\cdot\n=0 \; \; \mbox{for} \;\; y=0\; \mbox{and} \; x=0;  \qquad \u\cdot\n= \gamma e^{\pi x} \; \; \mbox{for} \;\; y=1;  \qquad \u\cdot\n= \gamma e^{\pi y} \; \; \mbox{for} \;\; x=1.$
\end{enumerate}
{  Here, $\gamma$ is a parameter.}\\

{
To study the dependency of the convergence with the parameter $\alpha$, we consider $N=60$, $\beta=20$, $\gamma=20$, and for each $\alpha$, we stop the algorithm \eqref{V1hi} when the error $Err_L <1e^{-5}$. We consider that the algorithms don't converge if this condition is not reached after $10000$ iterations. \\
Tables \ref{tab111} and \ref{tab112} show, for $\u_h^0=\0$, the error $Err$ and the number of iterations $Nbr$ which describe the convergence of Algorithm \eqref{V1hi} with respect to $\alpha$ and for each example. We remark that Algorithm (\ref{V1hi}) always converges for Example \eqref{Ex1FA} and the best convergence is for $\alpha=100$, while it converges for $\alpha > 10$ for Example \eqref{Ex2FA} and the best convergence is for $\alpha=1000$.
\begin{table}[h!]

\begin{tabular}{|l|l|l|l|l|l|l|l|}
\hline
\bf $\alpha$  & \bf 0.001 & \bf .01 & \bf .1 & \bf 1 & \bf 10 & \bf 100 & \bf 1000  \\
\hline
 \bf Nbr  & 9562 & 9225 & 8513 & 3461 & 500 & 54 & 62  \\
\hline
 \bf Err   & -1.6137 & -1.6137 & -1.6137 & -1.6137 & -1.6137 & -1.6137 & -1.61367\\
\hline
\end{tabular}
\caption{Error $Err$ (in logarithmic scale) and number of iterations $Nbr$ for $\u_h^0=\0$ and for each $\alpha$ associated to Example \eqref{Ex1FA}  with Algorithm (\ref{V1hi}). ($\beta=20$, $\gamma=20$).}\label{tab111}
%
%
\begin{tabular}{|l|l|l|l|l|l|l|l|l|l|}
\hline
\bf $\alpha$ & \bf 0.001 & \bf .01 & \bf .1 & \bf 1 & \bf 10 & \bf 100 & \bf 1000 \\
\hline
 \bf Nbr  & >10000 & >10000 & >10000 & >10000 & 5002 & 550 & 51  \\
\hline
 \bf Err &  div & div & div & div & -0.568261 & -0.56826 & -0.568258  \\
\hline
\end{tabular}
\caption{Error $Err$ { (in logarithmic scale)} and number of iterations $Nbr$ for $\u_h^0=\0$ and for each $\alpha$ associated to Example \eqref{Ex2FA}  with Algorithm (\ref{V1hi}). ($\beta=20$, $\gamma=20$).}\label{tab112}
\end{table}

Tables \ref{tab2211} and \ref{tab2212} show, for $\u_h^0$ computed with the Darcy's problem, the error $Err$ and the number of iterations $Nbr$. We remark that Algorithm (\ref{V1hi}) always converges and the best convergence in both cases is for $\alpha=1000$. The main advantage in this case where $\u_h^0$ computed with the Darcy's problem is that the algorithms (\ref{V1hi}) converge for all considered values of $\alpha$ and the number of iterations $Nbr$ is less than that obtained with $\u_h^0=\0$. We can deduce that in these examples, the convergence condition $\alpha > C h^{-d}$ with the initial guess $\u_h^0=\0$ can be balanced and compensated by using the initial guess $\u_h^0$ computed with the Darcy problem.
\begin{table}[h!]

\begin{tabular}{|l|l|l|l|l|l|l|l|}
\hline
\bf $\alpha$  & \bf 0.001 & \bf .01 & \bf .1 & \bf 1 & \bf 10 & \bf 100 & \bf 1000  \\
\hline
 \bf Nbr  & 5046 & 4979 & 4216 & 1758 & 261 & 30 &  28 \\
\hline
 \bf Err   & -1.6137 & -1.6137 & -1.6137 & -1.6137 & -1.6137 & -1.6137 & -1.61368\\
\hline
\end{tabular}
\caption{Error $Err$ (in logarithmic scale) and number of iterations $Nbr$ for $\u_h^0$ solution of Darcy problem and for each $\alpha$ associated to Example \eqref{Ex1FA}  with Algorithm (\ref{V1hi}). ($\beta=20$, $\gamma=20$).}\label{tab2211}
%
%
\begin{tabular}{|l|l|l|l|l|l|l|l|}
\hline
\bf $\alpha$  & \bf 0.001 & \bf .01 & \bf .1 & \bf 1 & \bf 10 & \bf 100 & \bf 1000  \\
\hline
 \bf Nbr  & 2761 & 2730 & 2417 & 1250 & 321 & 59 &  13 \\
\hline
 \bf Err   & -0.568261 & -0.568261 & -0.568261 & -0.568261 & -0.568261 & -0.568261 & -0.568262\\
\hline
\end{tabular}
\caption{Error $Err$ (in logarithmic scale) and number of iterations $Nbr$ for $\u_h^0$ solution of Darcy problem and for each $\alpha$ associated to Example \eqref{Ex2FA}  with Algorithm (\ref{V1hi}). ($\beta=20$, $\gamma=20$).}\label{tab2212}
\end{table}
}

{  For further study, we take $N=60$, $\beta=10$, $\gamma=1$, and we consider the initial guess $\u_h^0$ computed with the Darcy problem.} Tables \ref{tab11} and \ref{tab12} show the error $Err$ and the number of iterations $Nbr$. {  We notice that the best convergence (in term of number of iterations) is not obtained for $\alpha=1000$ as in the previous results but it} is obtained for $\alpha=10$ for Example \eqref{Ex1FA} and for $\alpha=100$ for Example \eqref{Ex2FA}. }
\begin{table}[h!]
\begin{tabular}{|l|l|l|l|l|l|l|l|}
\hline
\bf $\alpha$  & \bf 0.001 & \bf .01 & \bf .1 & \bf 1 & \bf 10 & \bf 100 & \bf 1000  \\
\hline
 \bf Nbr  & 138 & 136 & 116 & 47 & 11 & 44 & 179  \\
\hline
 \bf Err   & -1.85606 & -1.85606 & -1.85606 & -1.85606 & -1.85606 & -1.85602 & -1.85567\\
\hline
\end{tabular}
\caption{Error $Err$ (in logarithmic scale) and number of iterations $Nbr$ for each $\alpha$ associated to Example \eqref{Ex1FA}  of Algorithm (\ref{V1hi}). {  ($\beta=10$, $\gamma=1$).}}\label{tab11}
%
%
\begin{tabular}{|l|l|l|l|l|l|l|l|l|l|}
\hline
\bf $\alpha$ & \bf 0.001 & \bf .01 & \bf .1 & \bf 1 & \bf 10 & \bf 100 & \bf 1000 \\
\hline
 \bf Nbr  & 738 & 727 & 620 & 272 & 47 & 22 & 79  \\
\hline
 \bf Err &  -1.15854 & -1.15854 & -1.15854 & -1.15854 & -1.15854 & -1.15873 & -1.15841  \\
\hline
\end{tabular}
\caption{Error $Err$ { (in logarithmic scale)} and number of iterations $Nbr$ for each $\alpha$ associated to Example \eqref{Ex2FA}  of Algorithm (\ref{V1hi}). {  ($\beta=10$, $\gamma=1$).}}\label{tab12}
\end{table}

\begin{figure}[h!]
\begin{center}\vspace{-0cm}
   \begin{tabular}{c}\hspace{-0.5cm}
\centering\includegraphics[height=6.5cm,width=8.cm]{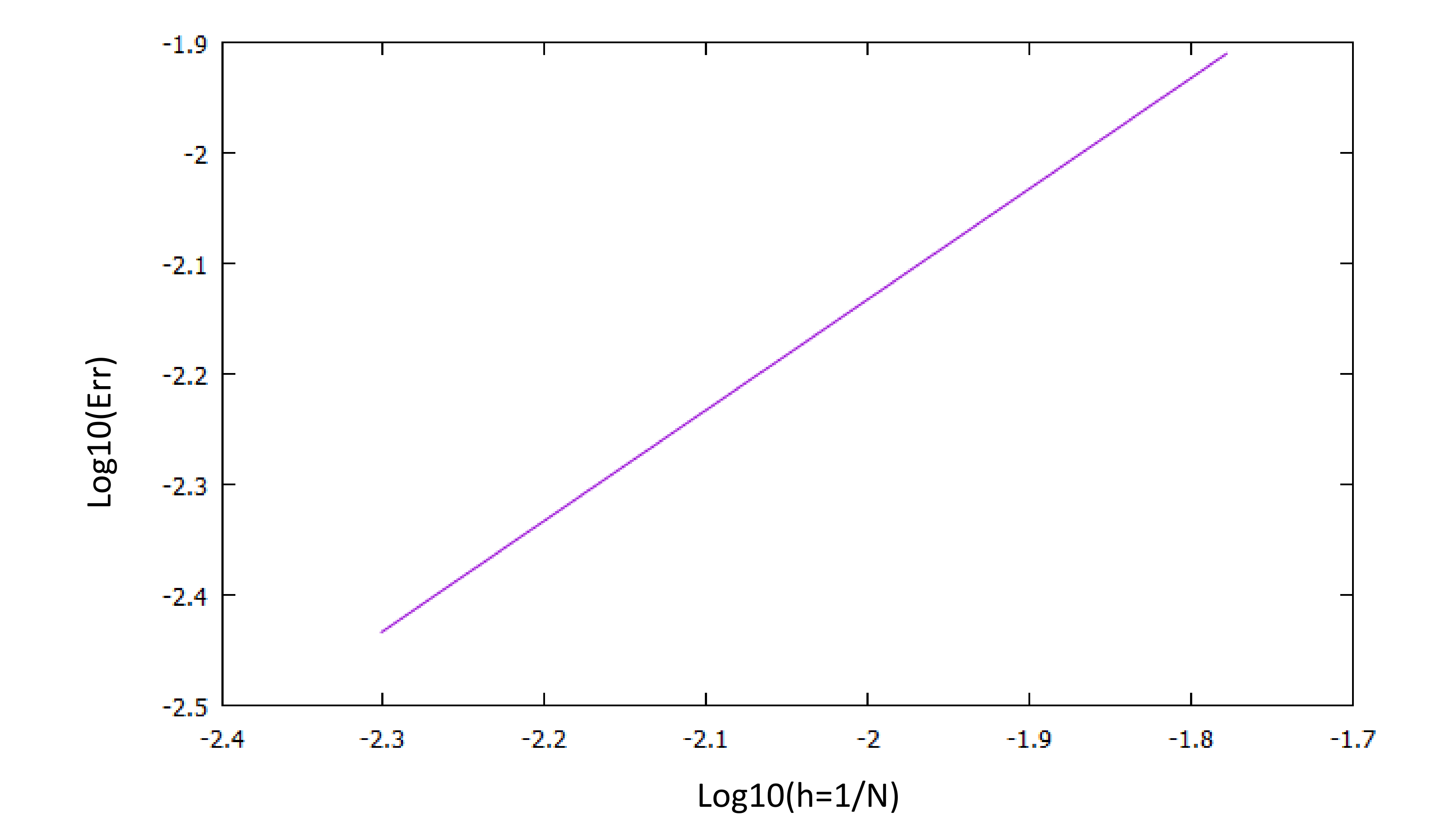}\hspace{0.cm}
\centering\includegraphics[height=6.5cm,width=8.cm]{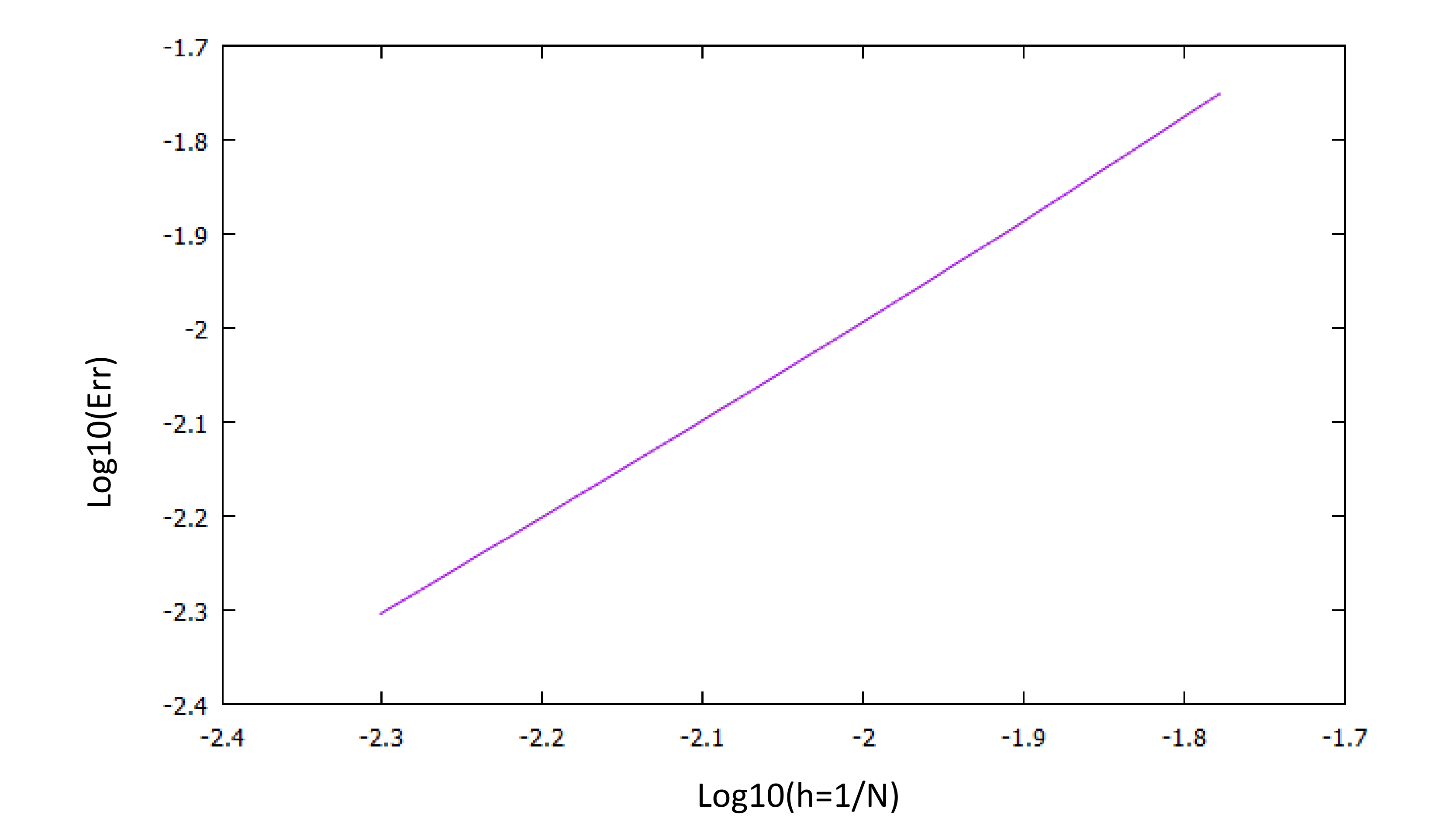}
\end{tabular}
\vspace{-.5cm}\caption{ the {\it  a priori} error for the algorithm \eqref{V1hi} with respect to $h=1/N$ :
left (first example) and right (second example). {  ($\beta=10$, $\gamma=1$).}}
\label{figV1ex12}
\end{center}
\end{figure}

Figures \ref{figV1ex12} shows, for $\alpha=10$ and $\gamma=1$ and $\beta=10$, in logarithmic scale the error $Err$ with respect to $h=\ds \frac{1}{N}, N=60, \dots , 200,$ for the algorithm \eqref{V1hi} (first example in the left and second example in the right). The slopes of the error lines are $1.001$ for the first example and $1.05$ for the second one.
%
%
%

{
\begin{rmq}\label{rmqq1}
The first discrete scheme \eqref{V1hi} corresponds to problem $(P1)$. We can list the  following comments:
\begin{enumerate}
\item The boundary condition must be on the normal component of the velocity  (boundary condition \eqref{E3}).
\item {  Tables \ref{tab111}, \ref{tab112}, \ref{tab2211}, \ref{tab2212}, \ref{tab11} and \ref{tab12}}  show that for the considered examples,
the number of iterations is relatively small when $\alpha$ is large. The disadvantage of this scheme is that we can not determine in advance the optimal value of $\alpha$ which leads to the minimal number of iterations. {  Furthermore, the initial guess $\u_h^0$ computed with Darcy problem gives better results then the initial guess $\u_h^0=\0$.}
\item The slopes of the curves presented in Figure \ref{figV1ex12} are close to the theoretical one (equal to $1$).
%
%
\end{enumerate}
\end{rmq}
}
{\rmq \label{remrkconvnum}
{  Theorems \ref{boundu1} and \ref{converg1}  show the convergence of Scheme \eqref{V1hi}  for $\alpha > C h^{-d}$ where $C$ can not be practically computable.
Tables \ref{tab111}, \ref{tab112}, \ref{tab2211} and \ref{tab2212} show that the initial guess computed with the Darcy problem gives better results. Furthermore, all the above tables allow us to deduce that for a given mesh (given $h$), we can consider in practice consider $\alpha \ge 10$ to test the convergence of the algorithm and we adjust the value of $\alpha$ in case of non-convergence.
}
}
\subsection{Second discrete scheme \eqref{V2hi}}
In this section we consider the scheme \eqref{V2hi} and we show numerical results. {  In fact, to compute the solution of the iterative problem \eqref{V2hi}, we use the penalty method (see \cite{Freefem}) which consists to solve the following problem:
\begin{equation*}
\left\{
\begin{array}{ll}
\medskip
\forall \v_h \in X_{p,h}, \quad \ds \int_\Omega \alpha (\u_h^{i+1} - \u_h^i) \cdot \v_h \, d\x +  \frac{\mu}{\rho} \int_\Omega K^{-1} \u^{i+1}_h \cdot \v_h \, d\x + \frac{\beta}{\rho} \int_\Omega |\u^i_h| \u^{i+1}_h \cdot \v_h \, d\x \\
\medskip
\hspace{5cm} - \ds \int_\Omega  p^{i+1}_h \, \div(\v_h) \, d\x = \ds \int_\Omega \f \cdot \v_h \, d\x,\\
\forall q_h\in M_{p,h}, \quad \ds  \int_\Omega  q_h  \div(\u^{i+1}_h) \, d\x  + \varepsilon \int_\Omega p_h^{i+1} \, q_h \, d\x = \int_\Omega b q_h \, d\x,
\end{array}
\right.
\end{equation*}
where $\varepsilon=10^{-8}$.\\
}

We propose the following two examples :
\begin{enumerate}
\item First example:
\begin{equation}\label{Ex1SA}
\left\{
\begin{array}{lcl}
\medskip
p(x,y) = 10\sin (\pi x) \sin (\pi y), &\qquad& \u(x, y) =  \gamma (\exp(x) \, \sin (\pi y),  \ds \frac{1}{\pi} \exp(x) \, \cos (\pi y))^T ,\\
b=0, &\qquad&  \f =	  \u +  \beta |\u| \u+ \nabla p,
\end{array}
\right.
\end{equation}
\item Second example:
\begin{equation}\label{Ex2SA}
\left\{
\begin{array}{lcl}
\medskip
p(x, y) = 10(x-x^2)(y-y^2), &\qquad& \u(x, y) =  \gamma (x \exp(\pi y), y \exp(\pi x))^T ,\\
b=\gamma(\exp(\pi x)+ \exp(\pi y)), &\qquad&  \f =	  \u +  \beta |\u| \u+ \nabla p,
\end{array}
\right.
\end{equation}
\end{enumerate}
{  Where $\gamma$ is a parameter.}\\

{  With this scheme \eqref{V2hi}, we repeat the same numerical tests as in the previous section for the scheme \eqref{V1hi}. We consider $N=60,  \beta=20$  and $\gamma=20$ and for each $\alpha$, we stop the algorithm \eqref{V2hi} when $Err_L <1e^{-5}$.\\
Tables \ref{tab211a} and \ref{tab212a} show for $\u_h^0=\0$ the error $Err$ which describes the convergence of Algorithm \eqref{V2hi} with respect to $\alpha$ and for each example. We remark that Algorithm (\ref{V2hi}) converges for Example \eqref{Ex1SA} for all the considered values of $\alpha$ and the best convergence is for $\alpha=1000$ while it converges for $\alpha > 10$ for Example \eqref{Ex2SA} and the best convergence is for $\alpha=1000$.
\begin{table}[h!]

\begin{tabular}{|l|l|l|l|l|l|l|l|}
\hline
\bf $\alpha$  & \bf 0.001 & \bf .01 & \bf .1 & \bf 1 & \bf 10 & \bf 100 & \bf 1000  \\
\hline
 \bf Nbr  & 8799 & 8666 & 7376 & 3020 & 487 & 61 & 45  \\
\hline
 \bf Err   & -1.4309 & -1.4309 & -1.43091 & -1.43091 & -1.43089 & -1.43089 & -1.43096\\
\hline
\end{tabular}
\caption{Error $Err$ (in logarithmic scale) and number of iterations $Nbr$ for $\u_h^0=\0$ and for each $\alpha$ associated to Example \eqref{Ex1SA}  of Algorithm (\ref{V2hi}). ($\beta=20$  and $\gamma=20$).}\label{tab211a}
%
%
\begin{tabular}{|l|l|l|l|l|l|l|l|l|l|}
\hline
\bf $\alpha$ & \bf 0.001 & \bf .01 & \bf .1 & \bf 1 & \bf 10 & \bf 100 & \bf 1000 \\
\hline
 \bf Nbr  & >10000 & >10000 & >10000 & >10000 & 3518 & 512 & 69  \\
\hline
 \bf Err &  div & div & div & div & -0.366231 & -0.366231 & -0.366229  \\
\hline
\end{tabular}
\caption{Error $Err$ { (in logarithmic scale)} and number of iterations $Nbr$ for $\u_h^0=\0$ and for each $\alpha$ associated to Example \eqref{Ex2SA}  of Algorithm (\ref{V2hi}). ($\beta=20$  and $\gamma=20$).}\label{tab212a}
\end{table}
}

{  Tables \ref{tab2211a} and \ref{tab2212a} show, for $\u_h^0$ computed with the Darcy's problem, the error $Err$ which describes the convergence of Algorithm \eqref{V2hi} with respect to $\alpha$ and for each example. In this case, we remark also that Algorithm (\ref{V2hi}) converges always and the best convergence is $\alpha=100$ for Example \eqref{Ex1FA} and $\alpha=1000$ for Example \eqref{Ex2FA}. The main advantage in this case where $\u_h^0$ computed with the Darcy's problem is that the number of iterations $Nbr$ is less than  that obtained with $\u_h^0=\0$.
\begin{table}[h!]

\begin{tabular}{|l|l|l|l|l|l|l|l|}
\hline
\bf $\alpha$  & \bf 0.001 & \bf .01 & \bf .1 & \bf 1 & \bf 10 & \bf 100 & \bf 1000  \\
\hline
 \bf Nbr  & 1218 & 1214 & 1038 & 463 & 104 & 18 &  20 \\
\hline
 \bf Err   & -1.4309 & -1.4309 & -1.4309 & -1.4309 & -1.43091 & -1.4309 & -1.43096\\
\hline
\end{tabular}
\caption{Error $Err$ (in logarithmic scale) and number of iterations $Nbr$ for $\u_h^0$ solution of Darcy problem and for each $\alpha$ associated to Example \eqref{Ex1FA}  of Algorithm (\ref{V1hi}). ($\beta=20$  and $\gamma=20$).}\label{tab2211a}
%
%
\begin{tabular}{|l|l|l|l|l|l|l|l|}
\hline
\bf $\alpha$  & \bf 0.001 & \bf .01 & \bf .1 & \bf 1 & \bf 10 & \bf 100 & \bf 1000  \\
\hline
 \bf Nbr  & 7618 & 7526 & 6781 & 3632 & 868 & 149 &  25 \\
\hline
 \bf Err   & -0.366232 & -0.366232 & -0.366229 & -0.366232 & -0.366232 & -0.366228 & -0.366228\\
\hline
\end{tabular}
\caption{Error $Err$ (in logarithmic scale) and number of iterations $Nbr$ for $\u_h^0$ solution of Darcy problem and for each $\alpha$ associated to Example \eqref{Ex2FA}  of Algorithm (\ref{V1hi}). ($\beta=20$  and $\gamma=20$).}\label{tab2212a}
\end{table}
}

{  For further study, we take $N=60$, $\beta=10$, $\gamma=1$, and we consider the initial guess computed with the Darcy problem.} Tables \ref{tab21} and \ref{tab22} show the error $Err$ and the number of iterations denoted by $Nbr$ (for the convergence) for each $\alpha$ and each example. {  We notice that the best convergence (in term of number of iterations) is  obtained for $\alpha=1000$ for Example \eqref{Ex1SA} and for $\alpha=100$ for Example \eqref{Ex2SA}. }

\begin{table}[h!]
\begin{tabular}{|l|l|l|l|l|l|l|l|}
\hline
\bf $\alpha$ &  \bf 0.001 & \bf .01 & \bf .1 & \bf 1 & \bf 10 & \bf 100 & \bf 1000  \\
\hline
 \bf Nbr  & 26 & 26 & 22 & 9 & 4 & 4 & 2  \\
\hline
 \bf Err  & -1.72478 & -1.72477 & -1.72477 & -1.72478 & -1.72477 & -1.72481 & -1.7252 \\
\hline
\end{tabular}
\caption{Error $Err$ (in logarithmic scale) and number of iterations $Nbr$ for each $\alpha$ associated to Example \eqref{Ex1SA}  of Algorithm (\ref{V2hi}). ($\beta=10$  and $\gamma=1$).}\label{tab21}
\end{table}
\begin{table}[h!]

\begin{tabular}{|l|l|l|l|l|l|l|l|}
\hline
\bf $\alpha$ &  \bf 0.001 & \bf .01 & \bf .1 & \bf 1 & \bf 10 & \bf 100 & \bf 1000  \\
\hline
 \bf Nbr & 490 & 481 & 376 & 181 & 40 & 17 & 59  \\
\hline
 \bf Err  & -0.96343 & -0.96345 & -0.96343 & -0.96347 & -0.96331 & -0.964411 & -0.963548  \\
\hline
\end{tabular}
\caption{Error $Err$ { (in logarithmic scale)} and number of {  iterations} $Nbr$ for each $\alpha$ associated to Example \eqref{Ex2SA}  of Algorithm (\ref{V2hi}).  ($\beta=10$  and $\gamma=1$).}\label{tab22}
\end{table}

\begin{figure}[h!]
\begin{center}\vspace{-0cm}
   \begin{tabular}{c}\hspace{-0.5cm}
\centering\includegraphics[height=6.5cm,width=8.cm]{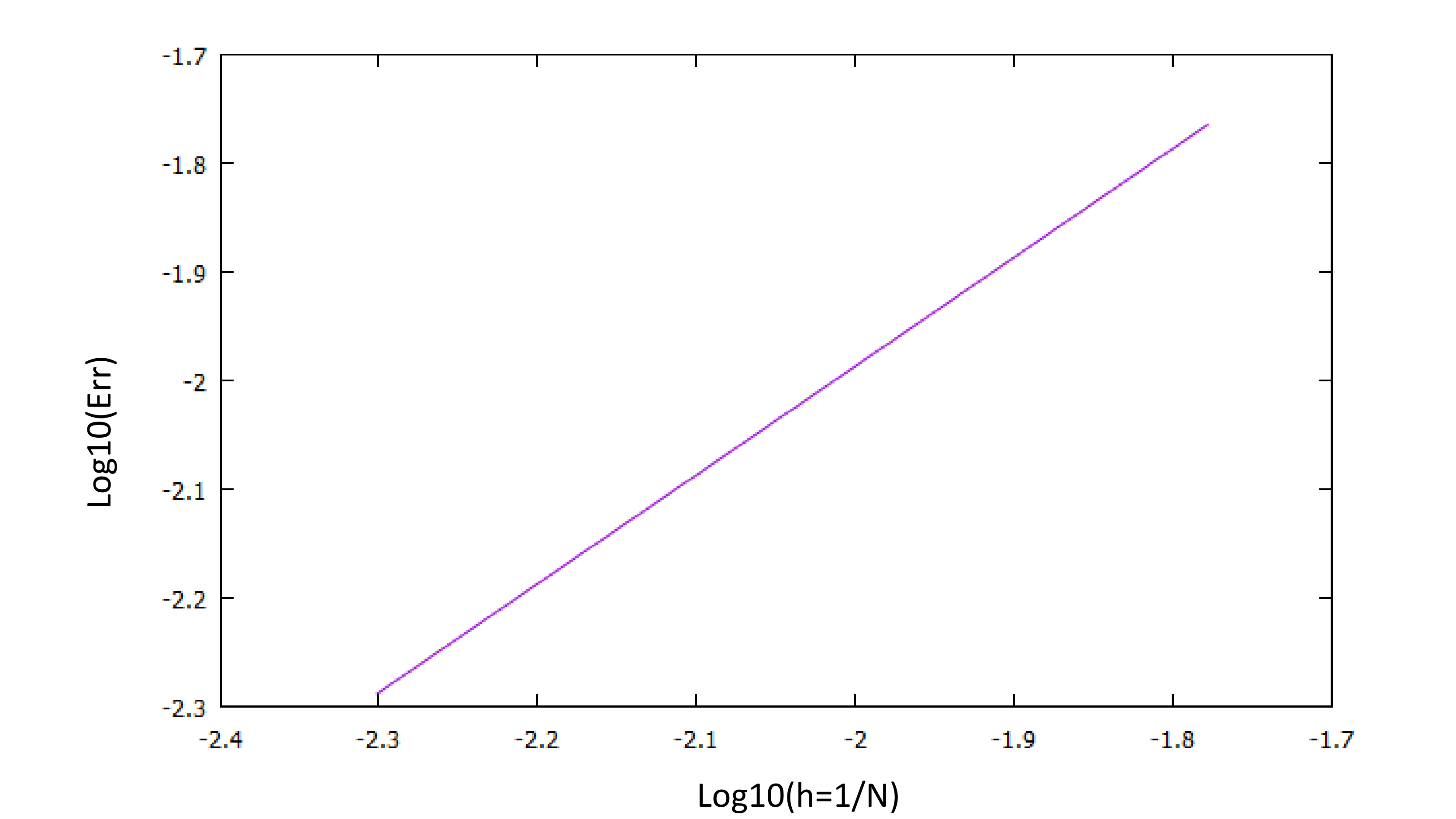}\hspace{0.cm}
\centering\includegraphics[height=6.5cm,width=8.cm]{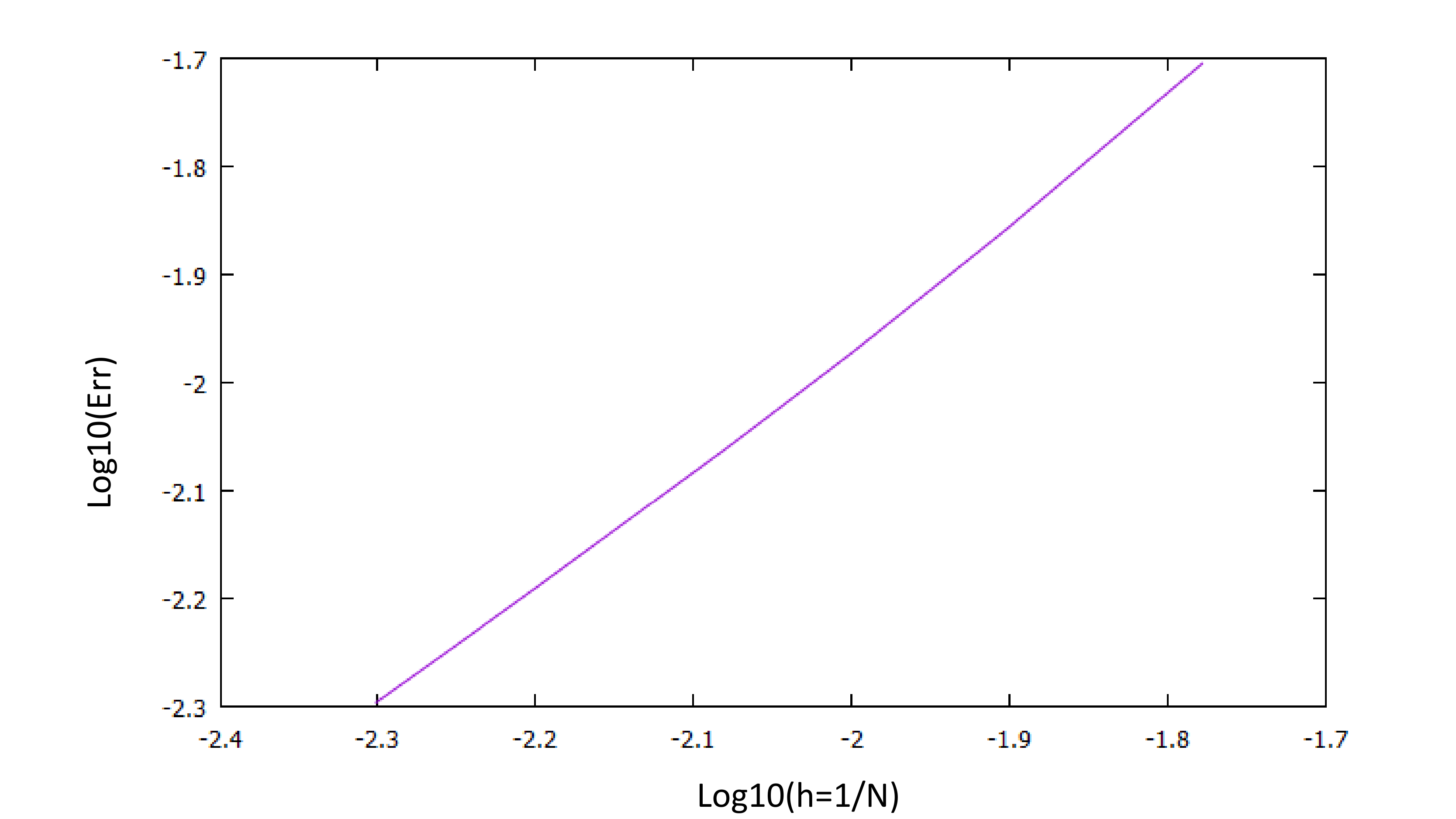}
\end{tabular}
\vspace{-.5cm}\caption{ A priori error for the algorithm (\ref{V2hi}) with respect to $h=1/N$ :
left (first example) and right (second example). ($\beta=10$  and $\gamma=1$).}
\label{figV2ex22}
\end{center}
\end{figure}

Figures \ref{figV2ex22} shows, for $\alpha=10$ and $\gamma=1$ and $\beta=10$, in logarithmic scale the error $Err$ with respect to $h=\ds \frac{1}{N}, N=60, \dots , 200,$ for the algorithm \eqref{V2hi} (first example in the left and second example in the right). The slopes of the error lines are $0.983$ for the first example and $1.10$ for the second one.\\

{  In order to study the dependence of Scheme \eqref{V2hi} with respect to the permeability tensor $K$,
we will consider the discontinuous tensor $K_\varepsilon$ (see \cite{exemp1K,exemp2K}):
\[
K(\x) = K_\varepsilon (\x)= \left\{
\begin{array}{rcl}
\varepsilon I && \quad \mbox{for all } x\in \Omega_0,\\
I && \quad \mbox{for all } x\in \Omega \backslash \Omega_0,
\end{array}
\right.
\]
where $I$ is the identity matrix, $\Omega_0 = [0.25, 0.5]\times  [0.25, 0.75] \subset \Omega$, and $\varepsilon = 10^{-6} \mbox{ or } 10^6$.\\
To study the dependency of the convergence with the parameter $\alpha$, we consider $\f=0$, $b=1$, $N = 60$ and $\beta=10$. Tables \ref{VehKdisc2aa} (respectively \ref{VehKdisc2a}) shows the convergence and the number of iterations of Scheme  \eqref{V2hi} with respect to $\alpha$ for $\varepsilon=10^6$ and for $\u_h^0=\0$ (respectively for $\u_h^0$ computed with Darcy's problem). We remark that the algorithm, with $\u_h^0=\0$, does not converge for $\alpha=.0001$ and converges for the other considered values of $\alpha$, while with $\u_h^0$ solution of Darcy problem, it converges always and the number of iterations is smaller for all the considered values of $\alpha$.
Furthermore, Tables \ref{VehKdisc2bb} (respectively \ref{VehKdisc2b}) shows the convergence and the number of iterations of Scheme \eqref{V2hi} with respect to $\alpha$ for $\varepsilon=10^{-6}$ and for $\u_h^0=\0$ (respectively for $\u_h^0$ computed with Darcy's problem). In this case, we remark that Algorithm \eqref{V2hi} converges always and the number of iterations is the same for the two cases of the initial guess $\u_h^0$.

\begin{table}[h!]

\begin{tabular}{|l|l|l|l|l|l|l|l|l|l|}
\hline
\bf $\alpha$ & \bf 0.001 & \bf .01 & \bf .1 & \bf 1 & \bf 10 & \bf 100 & \bf 1000 \\
\hline
 \bf Nbr  & div & 503 & 47 & 16 & 87 & 486 & 2370  \\
\hline
\end{tabular}
\caption{$\varepsilon=10^6$ and $\u_h^0=\0$. Number of iterations $Nbr$ for each $\alpha$ of Algorithm (\ref{V2hi}).}\label{VehKdisc2aa}
\begin{tabular}{|l|l|l|l|l|l|l|l|l|l|}
\hline
\bf $\alpha$ & \bf 0.001 & \bf .01 & \bf .1 & \bf 1 & \bf 10 & \bf 100 & \bf 1000 \\
\hline
 \bf Nbr  & 660 & 272 & 34 & 8 & 37  & 183 & 712 \\
\hline
\end{tabular}
\caption{$\varepsilon=10^6$ and $\u_h^0$ solution of Darcy. Number of iterations $Nbr$ for each $\alpha$ of Algorithm (\ref{V2hi}).}\label{VehKdisc2a}
\end{table}
\begin{table}[h!]

\begin{tabular}{|l|l|l|l|l|l|l|l|l|l|}
\hline
\bf $\alpha$ & \bf 0.001 & \bf .01 & \bf .1 & \bf 1 & \bf 10 & \bf 100 & \bf 1000 \\
\hline
 \bf Nbr & 20 & 20 & 17 & 8 & 23 & 124 & 592 \\
\hline
\end{tabular}
\caption{$\varepsilon=10^{-6}$ and $\u_h^0=\0$. Number of iterations $Nbr$ for each $\alpha$ of Algorithm (\ref{V2hi}).}\label{VehKdisc2bb}
\begin{tabular}{|l|l|l|l|l|l|l|l|l|l|}
\hline
\bf $\alpha$ & \bf 0.001 & \bf .01 & \bf .1 & \bf 1 & \bf 10 & \bf 100 & \bf 1000 \\
\hline
 \bf Nbr & 19 & 19 & 16 & 7 & 22 & 123 & 591 \\
\hline
\end{tabular}
\caption{$\varepsilon=10^{-6}$ and $\u_h^0$ solution of Darcy. Number of iterations $Nbr$ for each $\alpha$ of Algorithm (\ref{V2hi}).}\label{VehKdisc2b}
\end{table}
}
%
%
%
%
%
{
\begin{rmq}
The second discrete scheme \eqref{V2hi} corresponds to problem $(P2)$ where the boundary condition is \eqref{E4} (Dirichlet for the pressure). In this case, the same remarks given in Remarks \ref{rmqq1} and \ref{remrkconvnum} are also valid.
\end{rmq}
}
{\conc In this paper, we treat the Darcy-Forchheimer problem completed with two kinds of boundary conditions. We introduce two numerical schemes and show the corresponding convergences. we end this work with numerical investigations of validation.\\
}

\noindent {\bf Data availability statement  } Data sharing is not applicable to this article.
%
%

%
\end{document}